\documentclass[a4paper,12pt]{article}
\usepackage[bookmarks=true]{hyperref}
\usepackage{tikz}
\usepackage{amsmath,amsthm}
\usepackage{amssymb}
\usepackage{mathrsfs}
\newtheorem{theorem}{Theorem}[section]
\newtheorem{lemma}[theorem]{Lemma}
\newtheorem{remark}[theorem]{Remark}
\newtheorem{Def}[theorem]{Definition}
\newtheorem{cor}[theorem]{Corollary}
\newtheorem{pro}[theorem]{Proposition}

\numberwithin{equation}{section}

 \DeclareMathOperator{\R}{\mathbb
R} \textwidth=400pt

\def\mR{{\mathbb{R}}}

\def\mRT{{\mathbb{R}^2}}
\def\mRTh{{\mathbb{R}^3}}

\def\ourring{{\kappa \delta_{r_0,z_0} e_\theta}}

\def\defnv{{:=}}

\newcommand{\Real}{\mathbb R}

\newcommand{\be}{\begin{equation}}
\newcommand{\ee}{\end{equation}}

\newcommand{\rf}[1]{(\ref{#1})}
\newcommand{\pul}{{1\over 2}}

\newcommand{\om}{\omega}
\newcommand{\vf}{\varphi}

\newcommand{\curl}{{\rm curl\,}}
\newcommand{\intrt}{\int_{\Real^3}}

\def\bc{\begin{center}}       \def\ec{\end{center}}
\def\ba{\begin{array}}        \def\ea{\end{array}}
\def\be{\begin{equation}}     \def\ee{\end{equation}}
\def\bea{\begin{eqnarray}}    \def\eea{\end{eqnarray}}
\def\beaa{\begin{eqnarray*}}  \def\eeaa{\end{eqnarray*}}

\def\ben{\begin{enumerate}}   \def\een{\end{enumerate}}

         \def\om{\omega}
              
\def\oo{\infty}              
\def\la{\lambda}

\def\e{\varepsilon}

\def\lb{\label}              \def\x#1{(\ref{#1})}
\def\q{\quad}                \def\qq{\qquad}

\def\Proof{\noindent{\bf Proof}\, }
\def\qed{\hfill $\Box$ \smallskip}

\newcommand{\set}[1]{\left\{#1\right\}}
\newcommand{\norm}[1]{\left\Vert#1\right\Vert}
\newcommand{\abs}[1]{\left\vert#1\right\vert}

\begin{document}
\title{On the Cauchy problem for axi-symmetric vortex rings}
\author{Hao Feng\ \qq \ Vladim\'{i}r \v{S}ver\'{a}k}
\maketitle
\begin{abstract}
We consider the classical Cauchy problem for the 3d Navier-Stokes equation with the initial vorticity $\omega_0$ concentrated on a circle,
or more generally, a linear combination of such data for circles with common axis of symmetry.
We show that natural approximations of the problem obtained by smoothing the initial data satisfy good a-priori estimates
which enable us to conclude that the original problem with the singular initial distribution of vorticity has a solution.
We impose no restriction on the size of the initial data.
\end{abstract}

\section{Introduction}
Let us consider the classical Cauchy problem for the Navier-Stokes equation in $\Real^3\times(0,\infty)$:
\begin{eqnarray}
\left.\begin{array}{rcl}
u_t+{\rm div \,}(u \otimes u) +\nabla p -\nu\Delta u & = & 0 \\
{\rm div \,}u & = & 0
\end{array} \right\}& &  \quad \hbox{in $\Real^3\times(0,\infty)$}\,, \label{nse1} \\
u(\,\cdot\,,0)\,\,\, = \,\,\,u_0\,\,\,\, && \quad \hbox{in $\Real ^3$}\,.\label{nse2}
\end{eqnarray}
We will consider the initial data $u_0$ with vorticity $\om_0=\curl u_0$ which is supported on a circle.
In terms of the Geometric Measure Theory, $\om_0$ is a {\it \hbox{1-current}} of strength $\kappa$ supported on a smooth circle $\gamma$.
This means that for any smooth compactly supported test vector field (or, more precisely, 1-form) $\vf=(\vf_1,\vf_2,\vf_3)$ we can write
\be\lb{1}
\intrt \om_0\cdot \vf\,dx = \kappa\int_\gamma \vf_i(x)\,dx_i\,\,,
\ee
where the last integral is the classical curve integral (summation over the repeated indices is understood).
We will use the notation
\be\lb{2}
\om_0=\kappa\delta_\gamma
\ee
in this situation. The initial velocity field is recovered from $\om_0$ via the Biot-Savart law %\marginpar{{\tt please check the signs}}
\be\lb{3}
u_0(x)=-{1\over 4\pi} \intrt {(x-y)\wedge \om_0(y)\over \abs{x-y}^3}\,dy=-{1\over 4\pi}\int_\gamma  {(x-y)\wedge dy\over \abs{x-y}^3}\,\,.
\ee
We note that such $u_0$ has infinite kinetic energy:
\be\lb{4}
\intrt \pul|u_0|^2\,dx=+\infty\,\,,
\ee
due to the contributions from the immediate neighborhood of $\gamma$.
The initial datum of this form and its regularized variants are usually referred to as a {\it vortex ring}.
Their study goes back to Kelvin.
If $\gamma$ is the circle $(r_0\cos\theta,r_0\sin\theta,0)$ (with $-\pi\le\theta<\pi$) and $\kappa>0$,
we expect from Kelvin's calculations and the regularization due to the viscosity that at time $t$ the ring $\kappa\delta_\gamma$ will ``fatten" to thickness~$\sim \sqrt{\nu t}$ and will be moving up along the $z-$axis at speed roughly
\be\lb{5}
{\kappa\over 4\pi r_0} \log{a\over \sqrt{\nu t}}\,\,,
\ee
where $a$ is a suitable reference length.

Our goal here is to establish the existence of such a solution,
although we will not verify rigorously the detailed behavior suggested by Kelvin's calculations.
Our estimates will be less precise.
On the other hand, our method will be quite robust,
and can handle not only one vortex ring, but also a finite
or even continuous combination\footnote{with coefficients of the same sign} of such
as long as they have a common axis of symmetry.
The last condition is crucial, our method relies on the rotational symmetry of the situation.

It is instructive to compare our problem with the situation of parallel recti-linear vorticies. When the initial vorticity is supported on a line~$l$,
\be\lb{6}
\om_0=\kappa\delta_l\,,
\ee
the solution of the problem is given simply by the ``heat extension" of the initial data. When $l$ is the $x_3-$ axis, one has the text-book solution
\be\lb{7}
\om(x,t)=(0,0,\kappa \Gamma_2(x_1,x_2,\nu t))\,\,,
\ee
where $\Gamma_2(x_1,x_2,\nu t)={1\over 4\pi \nu t} e^{-{x_1^2+x_2^2\over 4\nu t}}$ is the 2d heat kernel. The non-linear term vanishes identically on these solutions.  Uniqueness is a subtle problem. The uniqueness has been proved in the class of the solutions of the form
\be\lb{7b}
u=(u_1(x_1,x_2,t), u_2(x_1,x_2,t),0)
\ee
(2d Navier-Stokes solutions), see \cite{GW2005,GGL2005}, but uniqueness among the 3d solutions seems to be open.

When the line $l$ is replaced by a collection of parallel lines $l_i$ and
\be\lb{8}
\om_0=\sum_i \kappa_i\delta_{l_i}
\ee
or possibly
\be\lb{9}
\om_0=\int \kappa_i\delta_{l_i}\,\,d\mu(i)\,\,,
\ee
where $\mu$ is a probability measure, one no longer has explicit solutions. The existence problem becomes more difficult
and was solved only in the 1980s in~\cite{Cottet1986,GMO}, see also \cite{BA1994,K1994}.
Uniqueness is again a subtle issue and is known  only in the class~\rf{7b} of 2d solutions, see~\cite{GG2005}.

Another class of existence results was obtained in \cite{GM1989} for small data, see also \cite{Taylor1992}.
In those papers the authors proved both existence and uniqueness (in suitable classes of functions) of the Cauchy problem~\rf{nse1},~\rf{nse2} for example in the case when the initial data $u_0$ is
\be\lb{10}
\om_0=\kappa\delta_\gamma\,,
\ee
where $\gamma$ is a smooth closed curve and $\kappa$ is sufficiently small (with the notion of smallness depending on $\gamma$). These results are proved by perturbation theory, and also follow from later works based on perturbation theory, such as~\cite{KT2001}.

Our main result in this paper is the following:

\begin{theorem} \lb{t1}
Let $\gamma$ be a circle, $\kappa\in\Real$ and $\om_0=\kappa \delta_\gamma$. Then the Cauchy problem ~\rf{nse1},~\rf{nse2} for the initial data $u_0$ given by $\om_0$ has a global solution which is smooth for $t>0$. The initial condition is satisfied in the following weak sense:
for any $\vf \in C^\oo_0 \bigl (\R^3; \R^3 \bigr )$,
\be \lb{e1}
\lim_{t \to 0}\int_{\R^3} \om(x,t) \cdot \vf(x)\,\, dx=\int_{\R^3} \om_0(x) \cdot \vf(x)\,\, dx\,\,,
\ee
where $\om=\curl u$ is the vorticity field.
\end{theorem}

\smallskip
{\bf Remarks}

1. Our method can be used to
show that the same results hold when $\om_0=\int \kappa(\gamma)\delta_\gamma\,\,d\mu(\gamma)$,
where $\mu$ is a probability measure supported on the set of the circles with a given axis of symmetry,
and $\kappa(\gamma)\ge 0$ is an integrable function with respect to $\mu$.

2. The sense in which the initial condition $u_0$ is assumed is somewhat weak, see (1.1).
A more precise analysis than ours is needed to determine optimal convergence of $\om(\,\cdot\,,t)\to\om_0$ as $t\to 0_+$.
%\footnote{{\tt Temporary note: We might also mention the following definition, closely related to (1.1).
%Let us take the vorticity equation
%\be\lb{pom1}
%\om_t+u\nabla \om-\om\nabla u - \dd \om = 0\,.
%\ee
%We multiply it by a smooth test-field $\vf(x,t)=(\vf_1(x,t),\vf_2(x,t),\vf_3(x,t))$
%which is supported in $B_R\times[0,T]$ for some $R>0,T>0$ (but does not have to vanish at $t=0$) and integrate by parts.
%We obtain
%\be\lb{pom2}
%-\intrt \om_0(x)\vf(x,0)\,dx+\int_0^T\intrt [-\om\vf_t-\om u\nabla\vf+u\om\nabla\vf -\om\dd \vf ]=0
%\ee
%One can check that this equation makes sense and it satisfied by our solution.
%We can then say that the initial condition $\om_0$ is attained weakly.
%}}

We now outline the main ideas involved in the proof.
By using the following transformation
\be
u(x,t) \mapsto \nu u(x, \nu t)\,, \hspace{10pt} p(x,t) \mapsto \nu^2 p(x, \nu t)\,\,,
\ee
we can change the first equation in \x{nse1} to
\be
u_t+{\rm div \,}(u \otimes u) +\nabla p -\Delta u  =  0\,\,.
\ee
Therefore, without loss of generality, we can assume $\nu=1$. Let us work with the vorticity equation
(obtained by taking the curl of the Navier-Stokes equations)
\be\lb{vort}
\om_t+u\nabla \om-\om\nabla u=\Delta \om\,\,,
\ee
which simplifies significantly for the axi-symmetric velocity fields with no swirl which we will be considering.
The precise definition is as follows.
\begin{Def} {\rm (Axi-symmetric vector field)}.\lb{d1}
A vector field $u$ in $\R^3$ is axi-symmetric if there is a coordinate frame in which
it can be written as
\be
u=u_r(r,z) e_r+u_\theta(r,z) e_\theta+u_z(r,z) e_z\,\,,
\ee
where
\be\lb{vectors}
e_r=({x_1/r},{x_2/r},0),\quad e_\theta=({-x_2/r}, {x_1/r}, 0),\,\quad e_z=(0,0,1)\,
\ee
and $(r,\theta,z)$ are the usual cylindrical coordinates associated with the frame.
The components $u_r$, $u_\theta$ and $u_z$ are independent of $\theta$.
The component $u_\theta$ is referred to as the swirl component of the vector field $u$ (in the given frame).
If $u_\theta$ vanishes, we say that $u$ has no swirl.
\end{Def}
It is easy to check that
the curl of an axi-symmetric vector field $u=u_r e_r+u_z e_z$ with no swirl is of the form
\be
\om=\hbox{curl}\ u=\left (u_{r,z}-u_{z,r} \right ) e_\theta\,\,,
\ee
which has only the $e_\theta$ component,
where $u_{r,z}$ denotes the partial derivative $\partial u_r/\partial z$, etc.
We will seek the solution of \x{vort} in the form $\om=\om_\theta(r,z,t) e_\theta$
and the velocity field in the form $u=u_r(r,z,t)e_r+u_z(r,z,t) e_z$.
The vorticity equation \x{vort} simplifies to
\be\lb{1-16}
\left({\omega_\theta\over r}\right)_t+
u \nabla
\left({
{\omega_\theta}\over r}\right)=\Delta \left({\omega_\theta\over r}\right)+\frac{2}{r}\left({
{\omega_\theta}\over r}\right)_{,r}.
\ee
The right hand side of \x{1-16} can be interpreted as the Laplacian in $\Real^5=\{(y_1,\dots,y_4,z)\}$ on functions which depend only on $r=\sqrt{y_1^2+\dots +y_4^2}$ and $z$.
Therefore the quantity $\omega_\theta\over r$ satisfies a maximum principle, see Lemma \ref{t15}.

There are three main ingredients of the proof:

\smallskip
1. Nash-type estimates for the quantity ${\omega_\theta\over r}$ based on equation~\rf{1-16} and the div-free nature of the field $u$.
These estimates give a good decay of $\norm{\om_\theta(t)\over r}_{L^\oo_x(\R^3)}$ in terms of $t^{-\alpha}$ for suitable $\alpha>0$,
even when the initial condition for $\om_\theta$ is a Dirac distribution, see \x{e72}.

\smallskip
2. The use of the conservation of the vorticity flux and momentum,
which are respectively the quantities $\int \omega_\theta(r,z)\,dr\,dz$ and $\int r^2\omega_\theta(r,z)\,dr\,dz$.

\smallskip
3. Weighted inequalities for axi-symmetric fields with no swirl, such as
\be\lb{18}
\norm{u}_{L^\oo_x(\R^3)}
\le C \norm{r \om_\theta}_{L^1_x(\R^3)}^{\frac{1}{4}} \norm{\frac{\om_\theta}{r}}^{\frac{1}{4}}_{L^1_x(\R^3)} \norm{\frac{\om_\theta}{r}}^{\frac{1}{2}}_{L^\oo_x(\R^3)}\,.
\ee
%\begin{tb}
%Assume that $\om$ is a vector field on $\R^3$ such that
%\[
%\norm{r \om}_{L^1_x(\R^3)}<\oo, \ \norm{\frac{\om}{r}}_{L^1_x(\R^3)}<\oo, \ \norm{\frac{\om}{r}}_{L^\oo_x(\R^3)}<\oo.
%\]
%Let $u$ be the vector field constructed from $\om$ via the Biot-Savart Law,
%\be \lb{Xs1bs}
%u(x)= -\frac{1}{4\pi} \int_{\mRTh} \frac{x-y}{\abs{x-y}^3} \times \om(y) dy.
%\ee
%Then, for any $\frac{3}{2} <q \le 6$, we have
%\be \lb{S1TBEB1}
%\norm{u}_{L^q_x(\R^3)} \le C(q) \norm{r \om}_{L^1_x(\R^3)}^{\frac{1}{2}}\norm{\frac{\om}{r}}_{L^1_x(\R^3)}^{\frac{1}{q}-\frac{1}{6}}\norm{\frac{ \om}{r}}_{L^\oo_x(\R^3)}^{\frac{2}{3}-\frac{1}{q}},
%\ee
%where $C(q)$ is an absolute constant independent of $\om$.
%Assume in addition that $\om$ has only angular component, $\om=\om_\theta e_\theta$,
%for some function $\om_\theta=\om_\theta(r,z)$,
%or equivalently, that the vector field $u$ in \x{Xs1bs} is axi-symmetric with no swirl. Then
%we have
%\be \lb{S1TBEB2}
%\begin{split}
%&\norm{u}_{L^\oo_x(\R^3)}
%\le N \norm{r \om}_{L^1_x(\R^3)}^{\frac{1}{4}} \norm{\frac{\om}{r}}^{\frac{1}{4}}_{L^1_x(\R^3)} \norm{\frac{\om}{r}}^{\frac{1}{2}}_{L^\oo_x(\R^3)},
%\end{split}
%\ee
%where $N$ is an absolute constant independent of $\om$ and moreover %, $\vec u$ is continuous on $\R^3$, satisfying
%\be \lb{Xs1dec}
%\lim_{\abs{x} \to \oo} \abs{u(x)}=0.
%\ee
%\end{tb}

\smallskip
Step 1 is achieved by applying of Nash's techniques \cite{N58}
for estimates of equations with div-free drift.
In our case they cannot quite be used directly,
due to the singular behavior of the coefficients of $\frac{2}{r}\left({
{\omega_\theta}\over r}\right)_{,r}$ near the $z-$axis which give extra terms in the Nash-type estimates.
Fortunately, the terms have a good sign, see the second term on line 6 in \x{e69} in the proof of Lemma \ref{t19}. %{\tt give reference to specific estimates}.
Inequality~\rf{18} seems to be of independent interest,
and it gives information about $u$ in terms of $\omega_\theta$, the quantity for which we have the most control.

Combining the results 1--3, we can then proceed along similar lines as \cite{GMO}.
The uniqueness of the solutions from the above theorem seems to be a difficult open problem.
We conjecture that it is possible to prove uniqueness in some natural  classes of axi-symmetric solutions without swirl,
but uniqueness in the class of all reasonable 3d vector fields may be much harder to prove and one might perhaps even have counter-examples.
We plan to consider these topics in a future work.

\section{Weighted inequalities} \lb{S2}

In this section, we present some weighted inequalities.
%Those inequalities stem out of our study on the main problem: the Navier-Stokes equations with
%the vortex ring $\kappa \delta_\gamma$ as initial vorticity, see Section \ref{S3}.
%But to get those inequalities, one do not need to assume that the vector field $u$ is a solution of Navier-Stokes equations.
%Thus we devote an independent section to introduce them. We are now in such a position:
We will have a-priori bounds on three quantities related to the vorticity:
$\norm{r \om}_{L^1_x(\R^3)}$, $\norm{\frac{\om}{r}}_{L^1_x(\R^3)}$, $\norm{\frac{\om}{r}}_{L^\oo_x(\R^3)}$,
and our aim is to obtain further estimates on the velocity $u$ from these bounds.
The inequalities presented in this section will be sufficient for our purposes in this paper.
%As those inequalities are of independent interest, we devote
%a part of the paper to them.

\begin{pro} \lb{t2}
Let $f\colon \R^3 \to \R$ be such that $\norm{r f}_{L^1_x(\R^3)}$, $\norm{\frac{f}{r}}_{L^1_x(\R^3)}$ and $\norm{\frac{f}{r}}_{L^\oo_x(\R^3)}$ are finite,
where $r=\sqrt{x_1^2+x_2^2}$. Then for every $1 \le p \le 2$, $f \in L^p_x(\R^3)$ and
\[
\norm{f}_{L^p_x(\R^3)} \le \norm{r f}_{L^1_x(\R^3)}^{\frac{1}{2}} \norm{\frac{f}{r}}_{L^1_x(\R^3)}^{\frac{1}{p}-\frac{1}{2}} \norm{\frac{f}{r}}_{L^\oo_x(\R^3)}^{1-\frac{1}{p}
}\,\,.
\]
\end{pro}
\Proof. We first prove the two cases of $p=1$ and $p=2$ and then use interpolation to prove the other cases. We can write
\[
\begin{split}
&\int_{\R^3} \abs{f} dx=\int_{\R^3} r^{\frac{1}{2}}\abs{f}^{\frac{1}{2}}
\frac{\abs{f}^{\frac{1}{2}}}{r^{\frac{1}{2}}} dx \\
&\le \biggl ( \int_{\R^3} r\abs{f} dx\biggr)^{\frac{1}{2}} \biggl (\int_{\R^3} \frac{\abs{f}}{r} dx \biggr)^{\frac{1}{2}}
= \norm{r f}_{L^1_x(\R^3)}^{\frac{1}{2}}\norm{\frac{f}{r}}_{L^1_x(\R^3)}^{\frac{1}{2}}\,,
\end{split}
\]
which proves the case $p=1$.

Next we consider
\[
\begin{split}
&\biggl (\int_{\R^3} \abs{f}^2 dx \biggr)^{\frac{1}{2}}=\biggl (\int_{\R^3} r\abs{f}\frac{\abs{f}}{r} dx \biggr)^{\frac{1}{2}} \\
&\le \biggl (\int_{\R^3} r\abs{f} \norm{\frac{f}{r}}_{L^\oo_x(\R^3)} dx \biggr)^{\frac{1}{2}} =\norm{r f}_{L^1_x(\R^3)}^{\frac{1}{2}} \norm{\frac{f}{r}}_{L^\oo_x(\R^3)}^{\frac{1}{2}}\,,
\end{split}
\]
which proves the case $p=2$.

Let $1<p<2$. We have
\[
\begin{split}
&\norm{f}_{L^p_x(\R^3)} \le \norm{f}_{L^1_x(\R^3)}^{\frac{2}{p}-1}\norm{f}_{L^2_x(\R^3)}^{2-\frac{2}{p}} \\
&\le \Biggl ( \norm{r f}_{L^1_x(\R^3)}^{\frac{1}{2}}\norm{\frac{f}{r}}_{L^1_x(\R^3)}^{\frac{1}{2}}\Biggr )^{\frac{2}{p}-1} \Biggl (\norm{r f}_{L^1_x(\R^3)}^{\frac{1}{2}} \norm{\frac{f}{r}}_{L^\oo_x(\R^3)}^{\frac{1}{2}}  \Biggr )^{2-\frac{2}{p}} \\
&=\norm{r f}_{L^1_x(\R^3)}^{\frac{1}{2}} \norm{\frac{f}{r}}_{L^1_x(\R^3)}^{\frac{1}{p}-\frac{1}{2}} \norm{\frac{f}{r}}_{L^\oo_x(\R^3)}^{1-\frac{1}{p}}.
\end{split}
\]
\qed

\begin{remark} Under the assumption of Proposition \ref{t2}, one can not control $\norm{f}_{L^p_x(\R^3)}$ for $p>2$.
It is not hard to exhibit counterexamples.
\end{remark}

\begin{cor} \lb{t3}
Assume that $\om$ is a vector field on $\R^3$ such that
\be \lb{e2}
\norm{r \om}_{L^1_x(\R^3)}<\oo, \ \norm{\frac{\om}{r}}_{L^1_x(\R^3)}<\oo, \ \norm{\frac{\om}{r}}_{L^\oo_x(\R^3)}<\oo\,\,.
\ee
Let $u$ be the vector field constructed from $\om$ via the Biot-Savart Law,
\be \lb{e3}
u(x)= -\frac{1}{4\pi} \int_{\mRTh} \frac{x-y}{\abs{x-y}^3} \times \om(y)\,\,dy\,\,.
\ee
Then for any $\frac{3}{2}<q \le 6$, $u \in {L^q_x(\R^3)}$ and
\be \lb{e4}
\norm{u}_{L^q_x(\R^3)} \lesssim \norm{r \om}_{L^1_x(\R^3)}^{\frac{1}{2}}\norm{\frac{\om}{r}}_{L^1_x(\R^3)}^{\frac{1}{q}-\frac{1}{6}}\norm{\frac{ \om}{r}}_{L^\oo_x(\R^3)}^{\frac{2}{3}-\frac{1}{q}}\,\,.
\ee
\end{cor}
\Proof. By Proposition \ref{t2} and \x{e2}, for any $1 \le p \le 2$, we have
\be \lb{e5}
\norm{\om}_{L^p_x(\R^3)} \le \norm{r \om}_{L^1_x(\R^3)}^{\frac{1}{2}} \norm{\frac{\om}{r}}_{L^1_x(\R^3)}^{\frac{1}{p}-\frac{1}{2}} \norm{\frac{ \om}{r}}_{L^\oo_x(\R^3)}^{1-\frac{1}{p}}\,\,.
\ee
Then by the classical Hardy-Littlewood-Sobolev inequality (see for instance \cite{Stein93,Torchinsky86}), one can get
\[
\norm{u}_{L^q_x(\R^3)} \lesssim \norm{\om}_{L^p_x(\R^3)}, \hspace{10pt} {\rm for} \ \
p \in (1,3) \ \ {\rm and } \ \ \frac{1}{q}=\frac{1}{p}-\frac{1}{3}\,\,,
\]
which, combining with \x{e5}, implies \x{e4}. \qed

\begin{remark} By interpolation, the a-priori bounds \x{e2} imply
\[
\norm{\frac{\om}{r}}_{L^p_x(\R^3)}<\oo, \hspace{20pt} {\rm for \ all} \ 1<p<\oo\,\,.
\]
What can we say about the full gradient $\nabla u$ from the above bounds and \x{e2}?
This question is related to the theory of singular integral operators with weights. Here we will only consider this question for vector fields which are axi-symmetric.
%Maybe we can not hope much for general vector fields. Let's restrict the question to the class of axi-symmetric
%vector fields with no swirl.
%It seems that the standard theory
%does not work well on this problem (or, this problem does not fit well into the framework of the standard theory).
%If $p \ge 2$,  $r^{-p}$ is not an $A_p$ weight in $\R^3$.
%For definition of $A_p$ weight, see \cite{Stein93,Torchinsky86}.
%Even though $(1+r)^{-p}$ looks better than $r^{-p}$ at the $z$-axis, it is not an $A_p$ weight in $\R^3$ if $p > 2$.
\end{remark}

It is natural to ask whether we can control other $L^q_x(\R^3)$ norms of $u$ except $\frac{3}{2}<q \le 6$ under the assumptions
of Corollary \ref{t3}.
%, say the $L^\oo_x(\R^3)$ norm. To make this question ``correct'', we may need to add one more restriction
%to the vector field $u$ in \x{e3} by assuming furthermore that $u$ is axi-symmetric and swirl-free,
%because the a-priori bounds in \x{e2} actually come from our study of axi-symmetric vector fields with no swirl.
%We will give a positive answer to this question in the sequel. But our approach
%is not as simple as the proof of Corollary \ref{t3}. To give the reader some feeling, we start from a simple and heuristic
%weighted inequality \x{e6} (see below), where we strengthen the a-priori bounds \x{e2} to the full gradient.
The inequality \x{e6} below indicates what can be expected in this situation. We prove this inequality as a warm-up for the proof of our main inequality~\rf{18}.
%It might be interesting to compare \x{e6} with \x{S1TBEB2} in Theorem B in Introduction.
\begin{pro} \lb{t4}
Assume $f=f(x_1,x_2,z)=f \bigl (\sqrt{x_1^2+x_2^2},z \bigr )\colon \R^3 \to \R$ is smooth and
vanishes at infinity. Assume in addition that $\norm{r \nabla f}_{L^1_x(\R^3)}$, $\norm{\frac{\nabla f}{r}}_{L^1_x(\R^3)}$ and $\norm{\frac{\nabla f}{r}}_{L^\oo_x(\R^3)}$ are finite.
Then we have
\be \lb{e6}
\norm{f}_{L^\oo_x(\R^3)} \lesssim \norm{r\nabla f}_{L^1_x(\R^3)}^{\frac{1}{4}}\norm{\frac{\nabla f}{r}}_{L^1_x(\R^3)}^{\frac{1}{4}} \norm{\frac{\nabla f}{r}}_{L^\oo_x(\R^3)}^{\frac{1}{2}}.
\ee
\end{pro}
\Proof. Assume $\abs{f(r,z)}$ achieve its supremum
at $(r_0,z_0)$, that is,
\[
\norm{f}_{L^\oo_x(\R^3)}=\abs{f(r_0,z_0)}\,\,.
\]
By the boundedness of $\frac{\nabla f}{r}$, $\nabla f$ must vanish at $r=0$ (the $z$-axis). In particular,
$\nabla_z f=0$ along the $z$-axis. Thus, $f(0,z) \equiv 0$ by the assumption that $f$ vanishes at infinity.
Therefore, without loss of generality, we can assume $r_0>0$. By the fundamental theorem of calculus and the
H\"{o}lder's inequality
\[
\begin{split}
&\norm{f}_{L^\oo_x(\R^3)}=\abs{f(r_0,z_0)}=\abs{f(r_0,z_0)^2}^{\frac{1}{2}}=\abs{\int_{z_0}^\oo \partial_z f(r_0,z)^2\,\,dz}^{\frac{1}{2}} \\
&\lesssim \Bigl (\int_{z_0}^\oo \abs{f(r_0,z)} \abs{\partial_z f(r_0,z)}\,\,dz \Bigr )^{\frac{1}{2}} \\
&=\biggl (\int_{z_0}^\oo \abs{\int_{r_0}^\oo \partial_r f(r,z)\,\,dr} \abs{\partial_z f(r_0,z)}\,\, dz\biggr )^{\frac{1}{2}} \\
&\lesssim \biggl (\int_{z_0}^\oo \int_{r_0}^\oo \abs{\partial_r f(r,z)}\,\,dr \abs{\partial_z f(r_0,z)}\,\, dz\biggr )^{\frac{1}{2}} \\
&= \biggl (\int_{z_0}^\oo \int_{r_0}^\oo \Bigl ( r \abs{\partial_r f(r,z)}^{\frac{1}{2}}\abs{\partial_r f(r,z)}^{\frac{1}{2}}\frac{1}{r} \Bigr )\,
dr\,\abs{\partial_z f(r_0,z)}\,dz\biggr )^{\frac{1}{2}} \\
&\lesssim \biggl (\int_{z_0}^\oo \int_{r_0}^\oo \Bigl ( r \abs{\partial_r f(r,z)}^{\frac{1}{2}}\abs{\partial_r f(r,z)}^{\frac{1}{2}} \Bigr )\,
dr\,\frac{\abs{\partial_z f(r_0,z)}}{r_0}\,dz\biggr )^{\frac{1}{2}} \\
%&\lesssim \biggl (\int_{z_0}^\oo \int_{r_0}^\oo \Bigl ( r \abs{\partial_r f(r,z)}^{\frac{1}{2}}\abs{\partial_r f(r,z)}^{\frac{1}{2}} \Bigr )
%dr \frac{\abs{\partial_z f(r_0,z)}}{r_0} dz\biggr )^{\frac{1}{2}} \\
&\lesssim \biggl (\int_{z_0}^\oo \int_{r_0}^\oo r \abs{\partial_r f(r,z)}^{\frac{1}{2}}\abs{\partial_r f(r,z)}^{\frac{1}{2}}
\,dr\,dz\biggr )^{\frac{1}{2}} \biggl (\sup_z \frac{\abs{\partial_z f(r_0,z)}}{r_0} \biggr)^{\frac{1}{2}}  \\
%&\lesssim \biggl (\int_{-\oo}^\oo \int_{0}^\oo r \abs{\partial_r f(r,z)}^{\frac{1}{2}}\abs{\partial_r f(r,z)}^{\frac{1}{2}} dr dz\biggr )^{\frac{1}{2}} \biggl (\sup_{\R^3} \frac{\abs{\partial_z f(r,z)}}{r} \biggr)^{\frac{1}{2}}  \\
&\lesssim \biggl (\int_{-\oo}^\oo \int_{0}^\oo r^2 \abs{\partial_r f(r,z)}\,
dr\,dz\biggr )^{\frac{1}{4}}\biggl (\int_{-\oo}^\oo \int_{0}^\oo \abs{\partial_r f(r,z)}\,
dr\,dz\biggr )^{\frac{1}{4}} \biggl (\sup_{\R^3} \frac{\abs{\partial_z f(r,z)}}{r} \biggr)^{\frac{1}{2}}  \\
&\lesssim \norm{r\nabla f}_{L^1_x(\R^3)}^{\frac{1}{4}}\norm{\frac{\nabla f}{r}}_{L^1_x(\R^3)}^{\frac{1}{4}} \norm{\frac{\nabla f}{r}}_{L^\oo_x(\R^3)}^{\frac{1}{2}}.
\end{split}
\] \qed

%The a-priori bounds in Proposition \ref{t4} are given by the full %gradient $\nabla f$ (with weights), and
%so when estimating the function $f$, we express $f$ in terms of its %gradient.
%This idea is of course not new. One can be use it to prove the %Gagliardo-Nireberg-Sobolev inequality, see \cite{E1998}. But
%it is meaningful to our question in some sense. 
In light of \x{e6}, one might ask whether the following inequality is true.
\be \lb{e7}
\norm{u}_{L^\oo_x(\R^3)}
\lesssim \norm{r \om}_{L^1_x(\R^3)}^{\frac{1}{4}} \norm{\frac{\om}{r}}^{\frac{1}{4}}_{L^1_x(\R^3)} \norm{\frac{\om}{r}}^{\frac{1}{2}}_{L^\oo_x(\R^3)}\,\,.
\ee
%and further ask whether we can prove \x{e7} by using the Biot-Savart Law \x{e3}
%(note that both \x{e6} and \x{e7} are consistent under the scaling transformation).
%It is not clear to us how to prove \x{e7} by just using \x{e3}.
%Maybe \x{e7} is false for general vector fields with only the assumption \x{e2}.
We do not know whether~(\ref{e7}) is true for general vector fields, but we will show that it turns out to be true for the class of axi-symmetric vector fields with no swirl, which is enough for our purposes here.
%and to prove \x{e7}, we need to use their special structure.
%Besides the usual Biot-Savart Law \x{e3}, an axi-symmetric vector field with no swirl can be expressed in terms of its curl in another way,
We will use the {\it axi-symmetric Biot-Savart Law}. To introduce it,
we start from the so-called axi-symmetric stream function.

In cylindrical coordinates, the class of axi-symmetric vector fields with no swirl is in
the form $u=u_r(r,z)e_r+u_z(r,z) e_z$, see Definition \ref{d1}, and the divergence-free condition div$u=0$ turns out to be
\[
(ru_r)_{,r}+(ru_z)_{,z}=0\,\,,
\]
 which means that
\[
ru_r=-\psi_{,z}\,, \q \ ru_z=\psi_{,r}
\]
for a suitable function $\psi=\psi(r,z)$, called the axi-symmetric stream function, similar to the 2d situation.
Hence
\be \lb{e8}
u_r=-\frac{1}{r} \psi_{,z}, \ \ u_z=\frac{1}{r} \psi_{,r}\,\,.
\ee
It is easy to check that the curl of an axi-symmetric field $u$ with no swirl is in the form
\[
{\rm curl} u=\om_\theta e_\theta
\]
with $\om_\theta=u_{r,z}-u_{z,r}$. Therefore, we obtain
\[
L\psi \colon =-\frac{1}{r} \psi_{,rr}+\frac{1}{r^2} \psi_{,r}-\frac{1}{r} \psi_{,zz} =\om_\theta\,\,.
\]
The inverse operator $L^{-1}$ is given by
\be \lb{e9}
\psi(\bar r,\bar z)=\int_{-\oo}^\oo\int_0^\oo\frac{\bar r r }{4\pi} \int_{0}^{2\pi}
\frac{\cos \varphi\,\, d \varphi}{\Bigl [ r^2+\bar r^2-2\bar r r\cos \varphi +(z-\bar z)^2 \Bigr ]^{\frac{1}{2}}}\,\om_\theta(r,z)\,dr\,dz\,.
\ee
For the axi-symmetric stream function and the derivation of \x{e9}, we refer the readers to \cite{Sverak2011}.
We can express \x{e9} somewhat more explicitly as
\be \lb{e10}
\begin{split}
\psi(\bar r,\bar z)&=\int_{-\oo}^\oo\int_0^\oo\frac{\sqrt{\bar r r}}{2\pi} \int_{0}^{\pi}
\frac{\cos \varphi\,\,d \varphi}{\Bigl [2(1- \cos \varphi) + \frac{(r-\bar r)^2+(z-\bar z)^2}{\bar r r}\Bigr ]^{\frac{1}{2}}}\,\om_\theta(r,z)\,dr\,dz \\
&=\int_{-\oo}^\oo\int_0^\oo\frac{\sqrt{\bar r r}}{2\pi} F \biggl (\frac{(r-\bar r)^2+(z-\bar z)^2}{\bar r r} \biggr)\,\om_\theta(r,z)\,dr\,dz\,,
\end{split}
\ee
where the function $F:(0,\oo) \to \R$ is defined by
\be \lb{e11}
F(s)\colon=\int_{0}^{\pi}
\frac{\cos \varphi\,\,d \varphi}{\bigl [2(1- \cos \varphi) + s\bigr ]^{\frac{1}{2}}}\,\,.
\ee
Let
\be \lb{e12}
G(\bar r, \bar z, r, z)=\frac{\sqrt{\bar r r}}{2\pi} F \biggl (\frac{(r-\bar r)^2+(z-\bar z)^2}{\bar r r} \biggr)\,\,.
\ee
Then
\[
\psi(\bar r,\bar z)=\int_{-\oo}^\oo\int_0^\oo G(\bar r, \bar z, r, z)\,\om_\theta(r,z)\,dr\,dz\,\,.
\]
By \x{e8} and \x{e12}, we get
\be \lb{e13}
\begin{split}
u_r(\bar r, \bar z)&=\int_{-\oo}^\oo\int_0^\oo \biggl [- \frac{1}{\bar r} \frac{\partial G}{\partial \bar z}(\bar r, \bar z, r, z) \biggr ]  \om_\theta(r,z)\,dr\,dz \\
%&=\int_{-\oo}^\oo\int_0^\oo \mathscr{R}(\bar r, \bar z, r,z)\om_\theta(r,z) dr dz
&=\int_{-\oo}^\oo\int_0^\oo \frac{z-\bar z}{\pi \bar r^{\frac{3}{2}}\sqrt{r}}  F^\prime \biggl ( \frac{(r-\bar r)^2+(z-\bar z)^2}{\bar r r}\biggr) \om_\theta (r,z)\,dr\,dz\,, \\
\end{split}
\ee
\be \lb{e14}
\begin{split}
u_z(\bar r, \bar z)&=\int_{-\oo}^\oo\int_0^\oo \biggl [\frac{1}{\bar r} \frac{\partial G}{\partial \bar r}(\bar r, \bar z, r, z) \biggr ]\,\om_\theta(r,z)\,dr\,dz \\
&=\int_{-\oo}^\oo\int_0^\oo \mathscr{Z}(\bar r, \bar z, r, z)\,\om_\theta(r,z)\,dr\,dz\,\,,
\end{split}
\ee
where
\[
%\mathscr{R}(\bar r, \bar z, r,z)=- \frac{1}{\bar r} \frac{\partial G}{\partial \bar z}(\bar r, \bar z, r, z), \ \
\mathscr{Z}(\bar r, \bar z, r,z)=\frac{1}{\bar r} \frac{\partial G}{\partial \bar r}(\bar r, \bar z, r, z)\,\,.
\]
The formulae \x{e13} and \x{e14}, representing the relations between $u_r, u_z$ and $\om_\theta$, represent the {\it axi-symmetric Biot-Savart Law}.
We calculate the kernel $\mathscr{Z}$. Let $d^2=(r-\bar r)^2+(z-\bar z)^2$. Let $\xi=\xi(\bar r, \bar z, r, z)=\frac{d}{\sqrt{\bar r r}}$.
Then by \x{e12}, we have
\[
G(\bar r, \bar z, r, z)=\frac{d}{2\pi \xi}F(\xi^2)=\frac{d}{2\pi} H(\xi)\,\,,
\]
where $H(t)=\frac{F(t^2)}{t}$. Direct calculation shows that
\be \lb{e15}
H^\prime(t)=2F^\prime(t^2)-\frac{F(t^2)}{t^2}\,, \ \  \frac{\partial \xi}{\partial \bar r}=\xi \Bigl (\frac{\bar r-r}{d^2}-\frac{1}{2 \bar r} \Bigr)\,\,,
\ee
\be \lb{e16}
\begin{split}
\mathscr{Z}=\frac{1}{\bar r}\frac{\partial G}{\partial \bar r}
&=\frac{1}{2\pi}\frac{\bar r-r}{\bar r^\frac{3}{2} r^\frac{1}{2}}\Bigl [ \frac{H(\xi)}{\xi}+H^\prime(\xi) \Bigr ]
-\frac{1}{4\pi} \xi^2 H^\prime(\xi) \frac{\sqrt{r}}{\bar r^\frac{3}{2}} \\
&=\frac{1}{\pi}\frac{\bar r-r}{\bar r^\frac{3}{2} r^\frac{1}{2}} F^\prime(\xi^2)
+\frac{1}{4\pi} \Bigl [ F(\xi^2)-2\xi^2 F^\prime(\xi^2) \Bigr ] \frac{\sqrt{r}}{\bar r^\frac{3}{2}}\,\,.
\end{split}
\ee
In the sequel, we are mainly interested in $\mathscr{Z}$ at $(\bar r, \bar z)=(1,0)$.
We write it down explicitly:
\be \lb{e17}
\begin{split}
\mathscr{Z}(1,0,r,z)
=&\frac{1-r}{\pi r^\frac{1}{2}} F^\prime \biggl ( \frac{(r-1)^2+z^2}{r}\biggr) \\
&+\frac{\sqrt{r}}{4\pi} \Biggl [ F\biggl ( \frac{(r-1)^2+z^2}{r}\biggr)-2\frac{(r-1)^2+z^2}{r} F^\prime \biggl ( \frac{(r-1)^2+z^2}{r}\biggr) \Biggr ]\,\,.
\end{split}
\ee
%\[ d=\sqrt{(r-\bar r)^2+(z-\bar z)^2} \]
At the first glance, comparing with the usual Biot-Savart Law \x{e3}, the axi-symmetric Biot-Savart Law \x{e13} and \x{e14} look more complicated and
have no advantages.
But \x{e13} and \x{e14} indeed capture some features of
axi-symmetric fields with no swirl. Although the function $F$ in \x{e11}
cannot be expressed in terms of elementary functions, it has nice asymptotic
properties near $s=0$ and $s=\oo$. By \x{e11}, it is obvious that
\be \lb{e18}
\abs{F(s)} \lesssim \Bigl (\frac{1}{s} \Bigr)^{\frac{1}{2}}.
\ee
However, $F$ actually has a
slower blow-up at $s=0$ and a faster decay at $s=\oo$ than \x{e18} as:
$\abs{F(s)} \lesssim \log \frac{1}{s}$ near
$s=0$ and $\abs{F(s)} \lesssim \Bigl (\frac{1}{s} \Bigr)^{\frac{3}{2}}$ near $s=\oo$.
We will use the following simple properties of $F$.
\begin{lemma}  \lb{t5}
For every non-negative integer $k$, the $k$th-derivative of $F$ satisfies
\be \lb{e19}
\abs{F^{(k)}(s)} \lesssim_k \ \frac{1}{s^{k+\frac{1}{2}}}\,\,,
\ee
for all $s \in (0,\oo)$.
\end{lemma}
\Proof. By \x{e11},
\[
\abs{F(s)} \lesssim \int_{0}^{\pi}
\frac{ d \varphi}{s^{\frac{1}{2}}} \lesssim \frac{1}{s^{\frac{1}{2}}}\,\,.
\]
Hence \x{e19} is true for the case of $k=0$. The first derivative of $F$ is
\[
F^\prime(s)=-\frac{1}{2} \int_{0}^{\pi}
\frac{\cos \varphi\,\, d \varphi}{\bigl [2(1- \cos \varphi) + s\bigr ]^{\frac{3}{2}}}\,\,.
\]
Therefore,
\[
\abs{F^\prime(s)} \lesssim \int_{0}^{\pi}
\frac{ d \varphi}{s^{\frac{3}{2}}} \lesssim \frac{1}{s^{\frac{3}{2}}}\,\,.
\]
Hence the case of $k=1$ is also true. The remaining cases can be proved similarly. \qed

\begin{lemma} \lb{t6}
There exists an absolute constant $0<\e_0<1$ such that for all $s \in (0,\e_0)$, the $k$th-derivative
of $F$ satisfies
\be \lb{e20}
\begin{split}
&\abs{F(s)} \lesssim \ \log \frac{1}{s} \lesssim_{\tau} \frac{1}{s^\tau}\,, \hspace{10pt} {\rm for \ every \ \tau>0, \ if} \ k=0\,, \\
&\abs{F^{(k)}(s)} \lesssim_{k} \ \frac{1}{s^k}\,, \hspace{95pt} {\rm if} \ 0<k \in \mathbb{N}\,.
\end{split}
\ee
\end{lemma}

\Proof. $F(s)$ has the following expansion near $s=0$, see for instance \cite{Sverak2011}
\[
F(s)=\Bigl (\log \frac{1}{s} \Bigr ) (a_0+a_1s+a_2s^2+\cdots)+(b_0+b_1s+b_2s^2+\cdots)\,\,,
\]
with $a_0=\frac{1}{2}$ and $b_0=\log 8-2$.
Hence
\[
F(s)=\frac{1}{2} \log \frac{1}{s}+\log 8-2+O\Bigl (s\log \frac{1}{s} \Bigr ), \hspace{20pt} s \to 0_+\,\,.
\]
The estimates \x{e20} follows easily from the above expansion.
\qed

\begin{lemma} \lb{t7}
There exists an absolute constant $N_0>1$ such that for every non-negative integer $k$, the $k$th-derivative of $F$ satisfies
\be \lb{e21}
\abs{F^{(k)}(s)} \lesssim_k \frac{1}{s^{k+\frac{3}{2}}}
\ee
for all $s \in (N_0, \oo)$.
\end{lemma}
\Proof. This is an easy calculation.
\qed

The estimates in Lemma \ref{t6} and Lemma \ref{t7} are local. But those restrictions
can be easily removed with the aid of Lemma \ref{t5}. As a consequence of Lemma \ref{t5}, Lemma \ref{t6} and Lemma \ref{t7}, we have
\begin{cor} \lb{t8}
For every non-negative integer $k$, the $k$th-derivative of $F$ satisfies
\[
\abs{F(s)} \lesssim_\tau \ \min \biggl (\Bigl ( \frac{1}{s}\Bigr )^\tau, \Bigl (\frac{1}{s} \Bigr )^{\frac{1}{2}}, \Bigl (\frac{1}{s} \Bigr )^{\frac{3}{2}} \biggr ),  \hspace{10pt} {\rm for \ every \ 0<\tau<\frac{1}{2}, \ if} \ k=0,
\]
\[
\abs{F^{(k)}(s)} \lesssim_k \ \min \biggl ( \Bigl ( \frac{1}{s}\Bigr )^k,\Bigl (\frac{1}{s} \Bigr )^{k+\frac{1}{2}}, \Bigl (\frac{1}{s} \Bigr )^{k+\frac{3}{2}} \biggr ),
\hspace{20pt} {\rm if} \ 0<k \in \mathbb{N},
\]
for all $s \in (0,\oo)$.
\end{cor}
%The estimates in Corollary \ref{t8} are actually optimal. They are easier to be used than Lemma \ref{t6} and Lemma \ref{t7}
%because we don't need to worry about whether the argument of $F$ (or, $F^\prime$) lie in the correct interval.
%In the sequel, we sometimes need the full strength of Lemma \ref{t6} and Lemma \ref{t7} but we sometimes do not.

With the aid of Corollary \ref{t8}, controlling $L^\oo_x(\R^3)$ of $u$ via the a-priori bounds \x{e2} becomes tractable.
We need the following technical lemma.
\begin{lemma} \lb{t9}
Let $f: \R^2 \to \R$ be such that $\norm{f}_{L^1(\R^2)}<\oo$ and $\norm{f}_{L^\oo(\R^2)}<\oo$.
Let $K: \R^2 \to \R$ be such that $\abs{K(x)} \le \frac{C}{\abs{x-x_0}}$ for
some positive constant $C$, some point $x_0 \in \R^2$ and for all $x \in \R^2$.
Then
\[
\abs{\int_{\R^2} K(x) f(x)\,dx} \le 2\sqrt{2\pi}C \norm{f}_{L^1(\R^2)}^\frac{1}{2}\norm{f}_{L^\oo(\R^2)}^\frac{1}{2}.
\]
\end{lemma}
\Proof. For any $\rho>0$, we have
\[
\begin{split}
\abs{\int_{\R^2} K(x) f(x)\,dx}
&\le \int_{\abs{x-x_0} \le \rho} \frac{C}{\abs{x-x_0}}\abs{f(x)}dx+\int_{\abs{x-x_0} > \rho} \frac{C}{\abs{x-x_0}}\abs{f(x)}dx \\
&\le 2\pi C \rho \norm{f}_{L^\oo(\R^2)}+\frac{C}{\rho}\norm{f}_{L^1(\R^2)}\,\,.
\end{split}
\]
After minimizing the last term, we can get the desired result. \qed

Since an axi-symmetric vector field $u$ with no swirl is of the form $u=u_r(r,z) e_r+u_z(r,z) e_z$,
to estimate the $L^\oo_x(\R^3)$ norm of $u$, it is enough to estimate the $L^\oo$
norms of $u_r$ and $u_z$ over the $rz$-plane $\Omega:=\set{r\ge 0, z \in \R}$.
We will use the following simple identities.
\[
\norm{r \om}_{L^1_x(\R^3)}=2\pi \norm{r^2 \om_\theta}_{L^1(\Omega)}, \
\norm{\frac{\om}{r}}_{L^1_x(\R^3)}=2 \pi \norm{\om_\theta}_{L^1(\Omega)}, \
\norm{\frac{\om}{r}}_{L^\oo_x(\R^3)}=\norm{\frac{\om_\theta}{r}}_{L^\oo(\Omega)}.
\]
We first estimate the $r$-component $u_r$.
\begin{pro} \lb{t10}
Let $u_r$ be given by the formula \x{e13} with $\om_\theta$ satisfying
\[
\norm{r^2 \om_\theta}_{L^1(\Omega)}<\oo, \ \norm{\om_\theta}_{L^1(\Omega)}<\oo, \ \norm{\frac{\om_\theta}{r}}_{L^\oo(\Omega)}<\oo\,\,.
\]
Then
\be \lb{e22}
\norm{u_r}_{L^\oo(\Omega)}
\le \ C_1 \norm{r^2 \om_\theta}_{L^1(\Omega)}^{\frac{1}{4}} \norm{\om_\theta}^{\frac{1}{4}}_{L^1(\Omega)} \norm{\frac{\om_\theta}{r}}^{\frac{1}{2}}_{L^\oo(\Omega)}\,\,,
\ee
where $C_1$ is an absolute constant.
\end{pro}
\Proof. The estimate \x{e22} is invariant under the scaling and the translation in the $z$ variable
\[
u_r(r,z) \mapsto u_r(\la r, \la z+z_0), \ \ \om_\theta(r,z) \mapsto  \la \om_\theta(\la r,\la z+z_0)
\]
for every $\la>0$ and every $z_0 \in \R$, and therefore it is enough to prove
\be \lb{e23}
\abs{u_r(1,0)}
\lesssim \ \norm{r^2 \om_\theta}_{L^1(\Omega)}^{\frac{1}{4}} \norm{\om_\theta}^{\frac{1}{4}}_{L^1(\Omega)} \norm{\frac{\om_\theta}{r}}^{\frac{1}{2}}_{L^\oo(\Omega)}\,\,.
\ee
By \x{e13}
\be \lb{e24}
u_r(1,0)=\int_{-\oo}^\oo\int_0^\oo \frac{z}{\pi \sqrt{r}}F^\prime \biggl (\frac{(r-1)^2+z^2}{r} \biggr) \om_\theta(r,z)\,dr\,dz\,\,.
\ee
We split the right hand side of \x{e24} into two parts. One is on the region
\[
I_1=\set{\frac{1}{2} \le r \le 2,-1 \le z \le 1}
\]
and the other on the complement $I_2=\Omega \setminus I_1$.

On $I_1$, by Corollary \ref{t8} (using $\abs{F^\prime(s)} \lesssim \frac{1}{s}$),
the kernel of \x{e24} can be estimated as
\[
\abs{\frac{z}{\pi \sqrt{r}}F^\prime \biggl (\frac{(r-1)^2+z^2}{r} \biggr)} \lesssim
\frac{\abs{z}}{\sqrt{r}} \frac{r}{(r-1)^2+z^2} \lesssim \frac{1}{\sqrt{(r-1)^2+z^2}}=\frac{1}{\abs{(r,z)-(1,0)}}\,\,.
\]
Therefore, by Lemma \ref{t9} and the fact that $r \sim 1$ on $I_1$, we obtain
\be \lb{e25}
\begin{split}
&\abs{\iint_{I_1} \frac{z}{\pi \sqrt{r}}F^\prime \biggl (\frac{(r-1)^2+z^2}{r} \biggr)\,\om_\theta(r,z)\,dr\,dz} \\
=&\abs{\iint \frac{z}{\pi \sqrt{r}}F^\prime \biggl (\frac{(r-1)^2+z^2}{r} \biggr)\,\om_\theta(r,z)\,\chi_{I_1}\,dr\,dz} \\
\lesssim &\norm{\om_\theta}^{\frac{1}{2}}_{L^1(I_1)}\norm{\om_\theta}^{\frac{1}{2}}_{L^\oo(I_1)}
\lesssim \norm{r^2 \om_\theta}_{L^1(I_1)}^{\frac{1}{4}} \norm{\om_\theta}^{\frac{1}{4}}_{L^1(I_1)}
\norm{\frac{\om_\theta}{r}}^{\frac{1}{2}}_{L^\oo(I_1)}\,\,,
\end{split}
\ee
where $\chi_{I_1}$ is the characteristic function of $I_1$.

On $I_2$, by Corollary \ref{t8}, (using $\abs{F^\prime(s)} \lesssim \Bigl ( \frac{1}{s}\Bigr )^\frac{5}{2}$),
the kernel of \x{e24} can be estimated as
\[
\abs{\frac{z}{\pi \sqrt{r}}F^\prime \biggl (\frac{(r-1)^2+z^2}{r} \biggr)} \lesssim
\frac{\abs{z}}{\sqrt{r}} \biggl (\frac{r}{(r-1)^2+z^2} \biggr)^\frac{5}{2}
\lesssim \frac{1}{(r-1)^2+z^2}\,\,,
\]
which is square-integrable on $I_2$.
Therefore, noting that
$\abs{\om_\theta}=r^\frac{1}{2}\abs{\om_\theta}^\frac{1}{4}\abs{\om_\theta}^\frac{1}{4}\frac{\abs{\om_\theta}^\frac{1}{2}}{r^\frac{1}{2}}$,
by H\"{o}lder's inequality, we obtain
\be \lb{e26}
\abs{\iint_{I_2} \frac{z}{\pi \sqrt{r}}F^\prime \biggl (\frac{(r-1)^2+z^2}{r} \biggr)\,\om_\theta(r,z)\,dr\,dz}
\lesssim \norm{r^2 \om_\theta}_{L^1(I_2)}^{\frac{1}{4}} \norm{\om_\theta}^{\frac{1}{4}}_{L^1(I_2)}
\norm{\frac{\om_\theta}{r}}^{\frac{1}{2}}_{L^\oo(I_2)}\,\,.
\ee
Clearly, \x{e24}, \x{e25} and \x{e26} imply \x{e23}. The proposition is proved. \qed

To estimate $u_z$, we need the following technical lemma.
\begin{lemma} \lb{t11}
Assume that $\om_\theta$ is a function on $\Omega$ satisfying
\[
\norm{r^2 \om_\theta}_{L^1(\Omega)}<\oo, \ \norm{\om_\theta}_{L^1(\Omega)}<\oo, \ \norm{\frac{\om_\theta}{r}}_{L^\oo(\Omega)}<\oo\,\,.
\]
Then
\be \lb{e27}
\int_{z=-\oo}^\oo \int_{r=2}^\oo \abs{\om_\theta(r,z)} \frac{r^2}{[(r-1)^2+z^2]^\frac{3}{2}}\,dr\,dz
\lesssim \norm{r^2 \om_\theta}_{L^1(\Omega)}^{\frac{1}{4}} \norm{\om_\theta}^{\frac{1}{4}}_{L^1(\Omega)} \norm{\frac{\om_\theta}{r}}^{\frac{1}{2}}_{L^\oo(\Omega)}\,\,.
\ee
\end{lemma}
We remark that the integral domain $\Omega$ of the right hand side of \x{e27} can be replaced by $\set{r \ge 2}$,
where $\set{r \ge 2}$ is shorthand for the set $\set{r \ge 2, z \in \R}$. But \x{e27} is enough for our purpose.

\Proof. We can't use the H\"{o}lder's inequality directly to get \x{e27} because on the region $\set{r \ge \abs{z}}$,
the weight $\frac{r^2}{[(r-1)^2+z^2]^\frac{3}{2}} \sim \frac{1}{[(r-1)^2+z^2]^\frac{1}{2}}$, which is not
square-integrable on that region. We introduce some notations. Let $d^2=r^2+z^2$ and $f(r,z)=\frac{\om_\theta(r,z)}{r}$.
To prove \x{e27}, it is enough to show
\be \lb{e28}
\int_{z=-\oo}^\oo \int_{r=2}^\oo \abs{f} \frac{r^3}{d^3}\,\,dr\,dz
\lesssim \norm{r^3 f}_{L^1(\Omega)}^{\frac{1}{4}} \norm{rf}^{\frac{1}{4}}_{L^1(\Omega)} \norm{f}^{\frac{1}{2}}_{L^\oo(\Omega)}.
\ee
By Cauchy-Schwartz inequality, we have
\[
\norm{r^2 f}_{L^1(\Omega)} \le \norm{r^3 f}_{L^1(\Omega)}^{\frac{1}{2}} \norm{rf}^{\frac{1}{2}}_{L^1(\Omega)}.
\]
Therefore, to prove \x{e28}, it is enough to prove
\be \lb{e29}
\int_{z=-\oo}^\oo \int_{r=2}^\oo \abs{f} \frac{r^3}{d^3}\,dr\,dz
\lesssim \norm{r^2 f}_{L^1(\set{r \ge 2})}^\frac{1}{2} \norm{f}^{\frac{1}{2}}_{L^\oo(\set{r \ge 2})},
\ee
since $\set{r \ge 2} \subset \Omega$.
We may assume that $f$ is a function supported in $\set{r \ge 2}$ and vanishing elsewhere in $\Omega$, otherwise,
we can just replace $f$ by $f\chi_{\set{r \ge 2}}$. Under this assumption, it is enough to prove
\be \lb{e30}
\norm{f \frac{r^3}{d^3}}_{L^1(\Omega)}
\lesssim \norm{r^2 f}_{L^1(\Omega)}^\frac{1}{2} \norm{f}^{\frac{1}{2}}_{L^\oo(\Omega)}.
\ee
For $\la>0$, let $f_\la(r,z)=\la^2 f(\la r, \la z)$. Clearly, $f_\la$ is supported on $\set{r \ge \frac{2}{\la}}$.
It is easy to check that for every $\la>0$, we have
\[
\norm{f_\la \frac{r^3}{d^3}}_{L^1(\Omega)}=\norm{f \frac{r^3}{d^3}}_{L^1(\Omega)},
\norm{f_\la}_{L^\oo(\Omega)}=\la^2 \norm{f}_{L^\oo(\Omega)},
\norm{r^2 f_\la}_{L^1(\Omega)}=\la^{-2} \norm{r^2 f}_{L^1(\Omega)}\,\,.
\]
We find $\la_0>0$ so that $\norm{f_{\la_0}}_{L^\oo(\Omega)}=\norm{r^2 f_{\la_0}}_{L^1(\Omega)}$.
By calculation,
\[
\la_0=\Biggl (\frac{\norm{r^2 f}_{L^1(\Omega)}}{\norm{f}_{L^\oo(\Omega)}} \Biggr )^\frac{1}{4}\,\,.
\]
To prove \x{e30}, it is enough to prove
\be \lb{e31}
\norm{f_{\la_0} \frac{r^3}{d^3}}_{L^1(\Omega)} \lesssim \norm{r^2 f_{\la_0}}_{L^1(\Omega)}+\norm{f_{\la_0}}_{L^\oo(\Omega)}\,\,.
\ee
We distinguish two cases $0<\la_0 \le 1$ and $\la_0>1$.

{\bf Case 1.} $0<\la_0 \le 1$.

By definition, $f_{\la_0}$ is supported on $\set{r \ge \frac{2}{\la_0}}$, which lies in $\set{r \ge 1}$.
On the support of $f_{\la_0}$, it is clear that $\frac{r^3}{d^3} \le 1 \le r^2$ and hence \x{e31} is true.

{\bf Case 2.} $\la_0>1$.

In this case, we have
\[
\begin{split}
\norm{f_{\la_0} \frac{r^3}{d^3}}_{L^1(\Omega)}
&\lesssim \int_{-\oo}^\oo \int_2^\oo \abs{f_{\la_0}}\,dr\,dz+ \norm{f_{\la_0}}_{L^\oo(\Omega)}\int_{-\oo}^\oo \int_{\frac{2}{\la_0}}^2\,\frac{r^3}{d^3}\,dr\,dz \\
&\lesssim \norm{r^2 f_{\la_0}}_{L^1(\Omega)}+\norm{f_{\la_0}}_{L^\oo(\Omega)}\,\,.
\end{split}
\]
Therefore \x{e31} is true. The lemma is proved. \qed

We now estimate the $z$-component $u_z$. The work for $u_z$ is similar to that for $u_r$ in Proposition \ref{t10} but some part have to be treated differently.
\begin{pro} \lb{t12}
Let $u_z$ be given by the formula \x{e14} with $\om_\theta$ satisfying
\[
\norm{r^2 \om_\theta}_{L^1(\Omega)}<\oo, \ \norm{\om_\theta}_{L^1(\Omega)}<\oo, \ \norm{\frac{\om_\theta}{r}}_{L^\oo(\Omega)}<\oo\,\,.
\]
Then
\be \lb{e32}
\norm{u_z}_{L^\oo(\Omega)}
\le \ C_2\norm{r^2 \om_\theta}_{L^1(\Omega)}^{\frac{1}{4}} \norm{\om_\theta}^{\frac{1}{4}}_{L^1(\Omega)} \norm{\frac{\om_\theta}{r}}^{\frac{1}{2}}_{L^\oo(\Omega)}\,\,,
\ee
where $C_2$ is an absolute constant.
\end{pro}
\Proof. Since the estimate \x{e32} is invariant under the scaling and the translation in the $z$ variable,
it is enough to prove
\be \lb{e33}
\abs{u_z(1,0)}
\lesssim \ \norm{r^2 \om_\theta}_{L^1(\Omega)}^{\frac{1}{4}} \norm{\om_\theta}^{\frac{1}{4}}_{L^1(\Omega)} \norm{\frac{\om_\theta}{r}}^{\frac{1}{2}}_{L^\oo(\Omega)}\,\,.
\ee
By \x{e14},
\be \lb{e34}
u_z(1,0)=\int_{-\oo}^\oo\int_0^\oo \mathscr{Z}(1, 0, r,z)\,\om_\theta(r,z)\,dr\,dz\,\,,
\ee
where $\mathscr{Z}(1, 0, r,z)$ is given by \x{e17} as
\be \lb{e35}
\begin{split}
\mathscr{Z}(1,0,r,z)
=&\frac{1-r}{\pi r^\frac{1}{2}} F^\prime \biggl ( \frac{(r-1)^2+z^2}{r}\biggr) \\
&+\frac{\sqrt{r}}{4\pi} \Biggl [ F\biggl ( \frac{(r-1)^2+z^2}{r}\biggr)-2\frac{(r-1)^2+z^2}{r} F^\prime \biggl ( \frac{(r-1)^2+z^2}{r}\biggr) \Biggr ] \\
:=& \mathscr{Z}_1(r,z)+\mathscr{Z}_2(r,z)\,\,.
\end{split}
\ee
We split the right hand side of \x{e34} into two parts. One is on the region
\[
I_1=\set{\frac{1}{2} \le r \le 2,-1 \le z \le 1}
\]
and the other on the complement $I_2=\Omega \setminus I_1$.

On $I_1$, by Corollary \ref{t8}, $\mathscr{Z}_1$ can be estimated as (using $\abs{F^\prime(s)} \lesssim \frac{1}{s}$)
\[
\abs{\mathscr{Z}_1(r,z)} \lesssim \frac{\abs{1-r}}{r^\frac{1}{2}} \frac{r}{(r-1)^2+z^2} \lesssim \frac{1}{\abs{(r,z)-(1,0)}}
\]
and $\mathscr{Z}_2$ can be estimates as (using $\abs{F(s)} \lesssim \Bigl ( \frac{1}{s}\Bigr )^\frac{1}{2}$
and $\abs{F^\prime(s)} \lesssim \Bigl ( \frac{1}{s}\Bigr )^\frac{3}{2}$)
\[
\abs{\mathscr{Z}_2(r,z)} \lesssim \sqrt{r} \Biggl [ \biggl ( \frac{r}{(r-1)^2+z^2} \biggr )^\frac{1}{2}
+\frac{(r-1)^2+z^2}{r} \biggl ( \frac{r}{(r-1)^2+z^2} \biggr )^\frac{3}{2} \Biggr ]
\lesssim \frac{1}{\abs{(r,z)-(1,0)}}\,\,.
\]
Therefore, by Lemma \ref{t9} and the fact that $r \sim 1$ on $I_1$, we obtain
\be \lb{e36}
\begin{split}
&\abs{\iint_{I_1} \mathscr{Z}(1,0,r,z)\,\om_\theta(r,z)\,dr\,dz} \\
\lesssim &\norm{\om_\theta}^{\frac{1}{2}}_{L^1(I_1)}\norm{\om_\theta}^{\frac{1}{2}}_{L^\oo(I_1)}
\lesssim \norm{r^2 \om_\theta}_{L^1(I_1)}^{\frac{1}{4}} \norm{\om_\theta}^{\frac{1}{4}}_{L^1(I_1)}
\norm{\frac{\om_\theta}{r}}^{\frac{1}{2}}_{L^\oo(I_1)}\,\,.
\end{split}
\ee

On $I_2$, by Corollary \ref{t8}, $\mathscr{Z}_1$ can be estimated as
(using $\abs{F^\prime(s)} \lesssim \Bigl ( \frac{1}{s}\Bigr )^\frac{5}{2}$)
\[
\abs{\mathscr{Z}_1(r,z)} \lesssim
\frac{\abs{1-r}}{r^\frac{1}{2}} \biggl (\frac{r}{(r-1)^2+z^2} \biggr)^\frac{5}{2}
\lesssim \frac{1}{(r-1)^2+z^2}\,\,,
\]
which is square-integrable on $I_2$.
Therefore, by H\"{o}lder's inequality, we obtain
\be \lb{e37}
\abs{\iint_{I_2} \mathscr{Z}_1(r,z)\,\om_\theta(r,z)\,dr\,dz}
\lesssim \norm{r^2 \om_\theta}_{L^1(I_2)}^{\frac{1}{4}} \norm{\om_\theta}^{\frac{1}{4}}_{L^1(I_2)}
\norm{\frac{\om_\theta}{r}}^{\frac{1}{2}}_{L^\oo(I_2)}\,\,.
\ee
Unfortunately, the foregoing argument of $\mathscr{Z}_1$ does not work for $\mathscr{Z}_2$ because $\mathscr{Z}_2$ is not square-integrable on the region $I_2$.
By Corollary \ref{t8}, the best estimate for $\mathscr{Z}_2$ on $I_2$ is
(using $\abs{F(s)} \lesssim \Bigl ( \frac{1}{s}\Bigr )^\frac{3}{2}$ and $\abs{F^\prime(s)} \lesssim \Bigl ( \frac{1}{s}\Bigr )^\frac{5}{2}$)
\be \lb{e38}
\abs{\mathscr{Z}_2(r,z)} \lesssim \sqrt{r} \Biggl [ \biggl ( \frac{r}{(r-1)^2+z^2} \biggr )^\frac{3}{2}
+\frac{(r-1)^2+z^2}{r} \biggl ( \frac{r}{(r-1)^2+z^2} \biggr )^\frac{5}{2} \Biggr ] \sim \frac{r^2}{[(r-1)^2+z^2]^\frac{3}{2}}\,\,.
\ee
To overcome this difficulty, we split the region $I_2$ into two parts,
``good'' part $I_{21}:=I_2 \cap \set{r \le 2}$ and ``bad'' part $I_{22}:=I_2 \cap \set{r>2}=\set{r>2}$.
By \x{e38}, $\mathscr{Z}_2$ is clearly square-integrable on $I_{21}$ and therefore by H\"{o}lder's inequality, we obtain
\be \lb{e39}
\abs{\iint_{I_{21}} \mathscr{Z}_2(r,z)\,\om_\theta(r,z)\,dr\,dz}
\lesssim \norm{r^2 \om_\theta}_{L^1(I_{21})}^{\frac{1}{4}} \norm{\om_\theta}^{\frac{1}{4}}_{L^1(I_{21})}
\norm{\frac{\om_\theta}{r}}^{\frac{1}{2}}_{L^\oo(I_{21})}\,\,.
\ee
On the ``bad'' part $I_{22}$, by Lemma \ref{t11} and \x{e38}, we have
\be \lb{e40}
\begin{split}
\abs{\iint_{I_{22}} \mathscr{Z}_2(r,z)\,\om_\theta(r,z)\,dr\,dz}
&\lesssim \int_{z=-\oo}^\oo \int_{r=2}^\oo \abs{\om_\theta(r,z)} \frac{r^2}{[(r-1)^2+z^2]^\frac{3}{2}}\,dr\,dz \\
&\lesssim \norm{r^2 \om_\theta}_{L^1(\Omega)}^{\frac{1}{4}} \norm{\om_\theta}^{\frac{1}{4}}_{L^1(\Omega)} \norm{\frac{\om_\theta}{r}}^{\frac{1}{2}}_{L^\oo(\Omega)}\,\,.
\end{split}
\ee
Clearly, \x{e34}, \x{e35}, \x{e36}, \x{e37}, \x{e39} and \x{e40} imply \x{e33}. The proposition is proved.
\qed

The following proposition concerns the decay as $\abs{x} \to \oo$.
\begin{pro} \lb{t13} Let $u=u_r e_r+u_z e_z$
with $u_r$ given by \x{e13} and $u_z$ given by \x{e14} and with $\om_\theta$ satisfying
\[
\norm{r^2 \om_\theta}_{L^1(\Omega)}<\oo, \ \norm{\om_\theta}_{L^1(\Omega)}<\oo, \ \norm{\frac{\om_\theta}{r}}_{L^\oo(\Omega)}<\oo\,\,.
\]
Then for every $\e>0$, there exists a $R>0$ such that for every $x \in \R^3$ with $\abs{x}>R$,
we have
\[
\abs{u(x)} \le \frac{\norm{r^2 \om_\theta}_{L^1(\Omega)}^\frac{1}{2} \norm{\om_\theta}_{L^1(\Omega)}^\frac{1}{2}}{2(\abs{x}-R)^2}+\frac{\e}{2}\,\,.
\]
In particular, we have
\[
\lim_{\abs{x} \to \oo} \abs{u(x)}=0\,\,.
\]
\end{pro}
\Proof. We can assume
\[
\norm{r^2 \om_\theta}_{L^1(\Omega)}>0, \ \norm{\om_\theta}_{L^1(\Omega)}>0, \ \norm{\frac{\om_\theta}{r}}_{L^\oo(\Omega)}>0\,\,,
\]
otherwise, $u \equiv 0$ and the assertions are obviously true.
For any $\e>0$, we can find a $R>0$ so that $\om_1\colon=\om_\theta \chi_{\set{r^2+z^2 \ge R^2}}$ satisfies
\[
\norm{\om_1}_{L^1(\Omega)}<\frac{\e^4}{16(C_1^2+C_2^2)^2 \norm{r^2 \om_\theta}_{L^1(\Omega)} \norm{\frac{\om_\theta}{r}}_{L^\oo(\Omega)}^2}\,\,,
\]
where $C_1$ and $C_2$ are the constants from Proposition \ref{t10} and \ref{t12}.
Let $\om_2=\om_\theta-\om_1$. Let $u_1$ and $u_2$ be the vector fields constructed from $\om_1$ and $\om_2$ via
\x{e13} and \x{e14}, respectively. Clearly, $u=u_1+u_2$.
By Proposition \ref{t10} and \ref{t12}, we have
\be \lb{e41}
\norm{u_1}_{L^\oo_x(\R^3)}
\le \ \sqrt{C_1^2+C_2^2} \norm{r^2 \om_1}_{L^1(\Omega)}^{\frac{1}{4}} \norm{\om_1}^{\frac{1}{4}}_{L^1(\Omega)} \norm{\frac{\om_1}{r}}^{\frac{1}{2}}_{L^\oo(\Omega)}
\le \frac{\e}{2}\,\,.
\ee
We can also express $u_2$ in terms of $\om_2$ via the Biot-Savart Law in Cartesian coordinates
\[
u_2(x)= -\frac{1}{4\pi} \int_{\R^3} \frac{x-y}{\abs{x-y}^3} \times \om_2\,e_\theta\,dy\,\,.
\]
Since $\om_2$ is supported in the ball $B_R(0)$, for any $\abs{x}>R$, we have
\be \lb{e42}
\abs{u_2(x)} \le \frac{1}{4 \pi} \frac{\norm{\om_2}_{L^1_x(\R^3)}}{(\abs{x}-R)^2}
=\frac{1}{2} \frac{\norm{r\om_2}_{L^1(\Omega)}}{(\abs{x}-R)^2}
\le \frac{\norm{r^2 \om_\theta}_{L^1(\Omega)}^\frac{1}{2} \norm{\om_\theta}_{L^1(\Omega)}^\frac{1}{2}}{2(\abs{x}-R)^2}\,\,.
\ee
Clearly, \x{e41} and \x{e42} imply the first assertion. The second assertion follows immediately from the first one.
\qed
\begin{remark}
In the statement of Proposition \ref{t13}, the $R$ depends not only on the norms
\be \lb{e43}
\norm{r^2 \om_\theta}_{L^1(\Omega)}, \ \norm{\om_\theta}_{L^1(\Omega)}, \ \norm{\frac{\om_\theta}{r}}_{L^\oo(\Omega)}
\ee
but also on the distribution of $\om_\theta$. For example, let $\om_\theta(r,z)=\chi_{\set{1 \le r \le 2, \abs{z} \le 1}}$.
Let $\om_\theta^{z_0}(r,z)=\om_\theta(r,z-z_0)$. Let $u^{z_0}=u_r^{z_0} e_r+u_z^{z_0} e_z$
be the vector field constructed from $\om_\theta^{z_0}$ via \x{e13} and \x{e14}.
%with $u_r^{z_0}$ given by \x{e13} and $u_z^{z_0}$ given by \x{e14}.
Obviously, we have
\[
\norm{r^2 \om_\theta^{z_0}}_{L^1(\Omega)}=\norm{r^2 \om_\theta}_{L^1(\Omega)}, \
\norm{\om_\theta^{z_0}}_{L^1(\Omega)}=\norm{\om_\theta}_{L^1(\Omega)}, \
\norm{\frac{\om_\theta^{z_0}}{r}}_{L^\oo(\Omega)}=\norm{\frac{\om_\theta}{r}}_{L^\oo(\Omega)}\,\,,
\]
\[
u_r^{z_0}(r,z)=u_r(r,z-z_0), \ \ u_z^{z_0}(r,z)=u_z(r,z-z_0)\,\,,
\]
but $u$ and $\set{u^{z_0}}_{z_0 \in \R}$ do not have a uniform decay since
the profile of $u^{z_0}$ is just the translation of that of $u$ by $z_0$ in the $z$-direction.
Nevertheless, they have the uniform decay rate in the $r$-direction. Actually, we can prove the following result that
for any $0<\e<\frac{1}{2}$ and any $x \in \R^3$ with $r=\sqrt{x_1^2+x_2^2} \ge 1$,
\be \lb{e44}
\abs{u(x)} \le \frac{C}{r^{\frac{1}{2}-\e}}\,\,,
\ee
where the constant $C$ depends only on the size of the norms in \x{e43}, see \cite{F13}.
But it is not clear whether \x{e44} is optimal.
\end{remark}

\section{A-priori estimates} \lb{S3}
In this section, we present the a-priori estimates for natural approximate solutions obtained by regularizing
the initial data, before which, we introduce the notations used.
The superscript ``${(\e)}$'' indicates the quantity (scalar or vector or tensor-valued) is induced by regularized initial data.
Sometimes we use a function $f=f(r,z)$ defined on $[0,\oo) \times \mR$ as a function defined on $\mRTh$
in the following way:
\[
f(x_1,x_2,z)=f \biggl (\sqrt{x_1^2+x_2^2},z \biggr ), \hspace{10pt} {\rm for} \ \ (x_1,x_2,z) \in \mRTh\,\,.
\]
Let us get back to our problem. The initial vorticity is
\be \lb{e45}
\om_0=\kappa \delta_\gamma\,\,,
\ee
where $\kappa \in \R$ and $\gamma$ is a circle. Without loss of generality,
we assume that $\gamma$ is $(r_0 \cos \theta, r_0 \sin \theta, z_0)$ for some $r_0>0$, $z_0 \in \R$ and $-\pi \le \theta <\pi$.
Then \x{e45} is equivalent to, in the sense of distribution,
\be \lb{e46}
\om_0=\ourring\,\,,
\ee
where $\delta_{r_0,z_0}$ is the Dirac mass at $(r_0,z_0)$ in the $rz$-plane.
We will search a solution in the class of
axi-symmetric velocity fields with no swirl, which have the form
\be \lb{e47}
u=u_r(r,z,t)e_r+u_z(r,z,t) e_z\,\,.
\ee
The related vorticity fields have the form
\be \lb{e48}
\om=\om_\theta(r,z,t) e_\theta
\ee
with $\om_\theta=u_{r,z}-u_{z,r}$. Note that
a solution of the form \x{e48} is formally compatible to the initial condition \x{e46}.
The equation for $\om_\theta$ is
\be \lb{e49}
\partial_t \om_\theta+u_r \om_{\theta,r}+u_z \om_{\theta,z}-\frac{u_r}{r}\om_\theta
=\om_{\theta,rr}+\frac{1}{r} \om_{\theta,r}-\frac{1}{r^2} \om_\theta+ \om_{\theta,zz}\,\,,
\ee
which can also be written as:
\be \lb{e50}
\partial_t \om_\theta + u \cdot \nabla \om_\theta-\frac{u_r}{r}\om_\theta=\Delta \om_\theta -\frac{1}{r^2} \om_\theta\,\,,
\ee
where
$\Delta=\frac{\partial^2}{\partial r^2}+\frac{1}{r}\frac{\partial}{\partial r}+\frac{1}{r^2}\frac{\partial^2}{\partial \theta^2}+\frac{\partial^2}{\partial z^2}$
is the scalar Laplacian in $\mRTh$, expressed in the cylindrical coordinates.
$u \cdot \nabla \om_\theta=u_1 \om_{\theta,1}+u_2 \om_{\theta,2}+u_z \om_{\theta,z}$
is equal to $u_r \om_{\theta,r}+u_z \om_{\theta,z}$. In terms of $\om_\theta$,
the initial condition \x{e46} can be formulated as:
\be \lb{e51}
\om_\theta(r,z,0)=\kappa \delta_{r_0,z_0}\,\,.
\ee
But we will not use either \x{e49} or \x{e50} in our method
because these two equations have a vortex-stretching term $-\frac{u_r}{r}\om_\theta$.
%which is not easily handled in the energy method.
%However we can ``eliminate'' the streching term from equation \x{e49} by first dividing \x{e49} by $r$ and then using the
%Leibniz rule to combine the drift term $u_r \frac{\om_{\theta,r}}{r}$ and the stretching term $-\frac{u_r}{r^2}\om_\theta$
%into a drift term $u_r \bigl (\frac{\om_\theta}{r} \bigr )_{,r}$. After that, we get the equation for 
It is easier to work with the quantity $\eta = \om_\theta/r$, which satisfies
\be \lb{e52}
\eta_t+u_r\eta_{,r}+u_z\eta_{,z}=\eta_{,rr}+\frac{3}{r}\eta_{,r}+\eta_{,zz}\,\,,
\ee
or
\be \lb{e53}
\eta_t+u \cdot \nabla \eta=\Delta \eta+\frac{2}{r}\eta_{,r}\,\,.
\ee

\begin{remark}  \lb{t14}
For a smooth vector field $u$, the apparent singularity of $\eta=\om_\theta /r$ is only
an artifact of the coordinate choice. The quantity $\eta$ is actually a smooth function,
even across the $z$-axis, as long as $u$ is smooth, see \cite{LW09}.
\end{remark}
\subsection{Regularized initial data}
In terms of $\eta$, the initia data \x{e51} reads:
\be \lb{e54}
\eta_0(r,z) \defnv \eta(r,z,0)=\frac{\om_\theta(r,z,0)}{r}=\frac{\kappa \delta_{r_0,z_0}}{r}=\frac{\kappa}{r_0} \delta_{r_0,z_0}\,\,.
\ee
The last equality of \x{e54} holds in the sense of distribution.
If we take an arbitrary test function $\psi=\psi(r,z)$, then
\[
\Bigl (\frac{\kappa \delta_{r_0,z_0}}{r}, \psi \Bigr )=\Bigl (\kappa \delta_{r_0,z_0}, \frac{\psi}{r} \Bigr )
=\kappa \frac{\psi(r_0,z_0)}{r_0}=\Bigl (\frac{\kappa}{r_0} \delta_{r_0,z_0}, \psi \Bigr )\,\,.
\]
Let $\phi:\mRT \to \mR$ be the standard mollifier such that $\phi \in C^\oo_0(B_1(0))$, $\phi \ge 0$ and $\int_\mRT \phi(y)\,dy=1$.
And Let $\phi^{(\e)}(y_1,y_2) \defnv \e^{-2} \phi(\frac{y_1}{\e}, \frac{y_2}{\e})$. Here and in the sequel, we assume $0<\e<\frac{r_0}{2}$.
We define $\eta_0^{(\e)}$ by
\be \lb{e55}
\eta_0^{(\e)}(r,z)\defnv \bigl( \phi^{(\e)}*\eta_0 \bigr ) (r,z)=\frac{\kappa}{r_0}\e^{-2} \phi \Bigl ( \frac{r-r_0}{\e}, \frac{z-z_0}{\e} \Bigr )\,\,.
\ee
Clearly,  for every $0<\e<\frac{r_0}{2}$, $\eta_0^{(\e)}$ has a compact support which stays away from the $z$-axis at least $\frac{r_0}{2}$.
It is easy to check
\be \lb{e56}
\begin{split}
&\pi \abs{\kappa} \le \norm{\eta_0^{(\e)}}_{L^1_x} \le 3\pi\abs{\kappa}\,\,, \\%:=C_1, \\
&\frac{\pi}{4}\abs{\kappa}r_0^2 \le \frac{2\pi \abs{\kappa}}{r_0} (r_0-\e)^3 \le \norm{r^2\eta_0^{(\e)}}_{L^1_x} \le \frac{2\pi \abs{\kappa}}{r_0} (r_0+\e)^3 \le \frac{27\pi}{4}\abs{\kappa}r_0^2\,\,.
\end{split}
\ee
\begin{remark} Note that $\norm{\eta_0^{(\e)}}_{L^1_x} \sim \abs{\kappa}$ and $\norm{r^2\eta_0^{(\e)}}_{L^1_x} \sim \abs{\kappa}r_0^2$.
The bounds for $\norm{\eta_0^{(\e)}}_{L^1_x}$ depends only on the strength $\abs{\kappa}$ of the ring $\ourring$ but
the bounds for $\norm{r^2\eta_0^{(\e)}}_{L^1_x}$ depends on both the strength and $r_0$.
%The ring $\ourring$ has a compact support. We then have such a bound $M_1$.
Nevertheless, they are both independent of $\e$ and will serve the a-priori bounds.
The inequalities in \x{e56} are dimensionally consistent.
\end{remark}
Corresponding to $\eta_0^{(\e)}$, the initial vorticity field $\om_0^{(\e)}$
and velocity field $u_0^{(\e)}$ are
\be \lb{e57}
\om_0^{(\e)}\defnv r\,\eta_0^{(\e)}\,e_\theta \hspace{10pt} {\rm and} \hspace{10pt} u_0^{(\e)}(x)\defnv -\frac{1}{4\pi} \int_{\mRTh} \frac{x-y}{\abs{x-y}^3} \times  \om_0^{(\e)}(y)\,dy\,\,,
\ee
respectively and $\om_0^{(\e)}$ has compact support.
\subsection{Approximate solutions for regularized initial data}
Obviously the velocity $u_0^{(\e)}$ in \x{e57}
is axi-symmetric and swirl-free. And for each $\e$, $u_0^{(\e)} \in H^k_x(\R^3)$ for any $k \ge 0$ and satisfies
\be \lb{e58}
{\rm div\,}u_0^{(\e)}=0, \hspace{20pt}
{\rm curl\,}u_0^{(\e)}=\om_0^{(\e)}\,\,.
\ee
\begin{remark} We don't have a uniform bound for $H^k_x(\R^3)$ norms of $u_0^{(\e)}$,
not even for the $L^2_x(\R^3)$ norms of $u_0^{(\e)}$.
\end{remark}
Then by the result of \cite{L68,UY68,LMNP99},
there exists a unique global-in-time smooth solution $u^{(\e)}$ for 3d Navier-Stokes equations satisfying
the initial condition
\be \lb{e59}
u^{(\e)}(0)=u_0^{(\e)}\,\,.
\ee
And moreover $u^{(\e)}$ is axi-symmetric with no swirl, that is, in cylindrical coordinates,
\[
u^{(\e)}=u^{(\e)}_r(r,z,t) e_r + u^{(\e)}_z(r,z,t) e_z\,\,.
\]
We shall show that a subsequence of $\Bigl \{ u^{(\e)} \Bigr \}_{0<\e<\frac{r_0}{2}}$
converges to a smooth solution with the ring $\ourring$ as initial vorticity.
Corresponding to $u^{(\e)}$, the vorticity field $\om^{(\e)}$ and the scalar quantity $\eta^{(\e)}$ are
\be \lb{e60}
\om^{(\e)} = {\rm curl} \ u^{(\e)}=\Bigl ( u^{(\e)}_{r,z}-u^{(\e)}_{z,r} \Bigr ) e_\theta \hspace{15pt}
{\rm and}\hspace{15pt} \eta^{(\e)} = \frac{u^{(\e)}_{r,z}-u^{(\e)}_{z,r}}{r}\,\,,
\ee
respectively. As a result of \x{e57}, \x{e58}, \x{e59} and \x{e60},
$\om^{(\e)}$ and $\eta^{(\e)}$ satisfy the initial data in \x{e57}
\be \lb{e61}
\om^{(\e)}(0)=\om_0^{(\e)}, \
\eta^{(\e)}(0)=\eta_0^{(\e)}\,\,.
\ee
By \x{e53} and Remark \ref{t14}, $\eta^{(\e)}$ is a smooth solution of the following equation:
\be \lb{e62}
\eta^{(\e)}_t+u^{(\e)} \cdot \nabla \eta^{(\e)}=\Delta \eta^{(\e)}+\frac{2}{r}\eta^{(\e)}_{,r}, \hspace{10pt} {\rm in} \ \ \mRTh \times (0, \oo)\,\,.
\ee

\subsection{A-priori estimates for approximate solutions}
The following lemma says that $\eta^{(\e)}$ enjoys the strong maximum priciple, which is
crucial for our arguments of obtaining the a-priori estimates.
\begin{lemma} \lb{t15}
If $\kappa>0$ {\rm (}or, $<0${\rm )}, then $\eta^{(\e)}(r,z,t)>0$ {\rm (}or, $<0${\rm )} for any $r \ge 0$, $z \in \mR$ and $t>0$.
\end{lemma}
\Proof. We just prove the case of $\kappa>0$. The case of $\kappa<0$ can be proved similarly.
We can not apply the maximum principle directly to \x{e62} since
the coefficient of $\frac{2}{r}\eta^{(\e)}_{,r}$ is singular. Recalling that
the Laplacian of a radially symmetric function $v(r)$ defined on $\R^n$ is $\Delta v=v^{\prime\prime}(r)+\frac{n-1}{r}v^\prime(r)$,
the right hand side of \x{e62} can be appropriately interpreted as the Laplacian in $\R^5$ and
we can recast \x{e62} in $\mR^5 \times (0, \oo)$. To this end, we introduce some notations.
Define
\[
\begin{split}
&{\hat \eta}^{(\e)}(x_1,x_2,x_3,x_4,z,t) \defnv \eta^{(\e)} \Bigl ( \sqrt{x_1^2+x_2^2+x_3^2+x_4^2},z,t \Bigr )\,\,, \\
&{{\hat u}}^{(\e)}(x_1,x_2,x_3,x_4,z,t) \\
&\defnv u^{(\e)}_r\Bigl ( \sqrt{x_1^2+x_2^2+x_3^2+x_4^2},z,t \Bigr ) {{\hat e}}_r
+ u^{(\e)}_z\Bigl ( \sqrt{x_1^2+x_2^2+x_3^2+x_4^2},z,t \Bigr ) {{\hat e}}_z\,\,,
\end{split}
\]
where
\[
r=\sqrt{x_1^2+x_2^2+x_3^2+x_4^2}, \ \
{{\hat e}}_r=\Bigl (\frac{x_1}{r}, \frac{x_2}{r}, \frac{x_3}{r}, \frac{x_4}{r}, 0 \Bigr ), \ \
{{\hat e}}_z=\Bigl (0,0,0,0,1 \Bigr )\,\,.
\]
Then by \x{e55}, \x{e61} and \x{e62}, we have
\[
\left \{
\begin{array}{l}
{\hat \eta}^{(\e)}_t+{{\hat u}}^{(\e)} \cdot \nabla_5 {\hat \eta}^{(\e)} =\Delta_5 {\hat \eta}^{(\e)},
\hspace{10pt} {\rm in} \ \ \mR^5 \times (0, \oo)\,\,, \\
{\hat \eta}^{(\e)}(0) \ge 0, \ {\rm and} \ \not\equiv 0 \hspace{10pt} {\rm in} \ \ \mR^5\,\,,
\end{array}
\right.
 \]
where,
\[
\nabla_5=\Bigl (\frac{\partial}{\partial x_1}, \frac{\partial}{\partial x_2},
\frac{\partial}{\partial x_3}, \frac{\partial}{\partial x_4}, \frac{\partial}{\partial z} \Bigr), \ \
\Delta_5=\frac{\partial^2}{\partial x_1^2}+\frac{\partial^2}{\partial x_2^2}+
\frac{\partial^2}{\partial x_3^2}+\frac{\partial^2}{\partial x_4^2}+\frac{\partial^2}{\partial z^2}\,\,.
\]
By strong maximum principle, we get
\[
{\hat \eta}^{(\e)}>0, \hspace{10pt} {\rm in} \ \ \mR^5 \times (0, \oo)\,\,,
\]
which implies
\[
\eta^{(\e)}>0\,\,.
\]
Thus the lemma is proved. \qed \\

One of the important a-priori estimates is the conservation of momentum.
\begin{lemma} \lb{t16}{\rm (Conservation of momentum)}. %For the constant $M_1$ in \x{e56}, we have
For all $t \ge 0$, we have
\be \lb{e63}
%\frac{2\pi \abs{\al}}{r_0} (r_0-\e)^3 \le
\norm{r \om^{(\e)}(t)}_{L^1_x}=\norm{r \om^{(\e)}(0)}_{L^1_x} %\le \frac{2\pi \abs{\al}}{r_0} (r_0+\e)^3
%\le M_1
\le \frac{27\pi}{4}\abs{\kappa}r_0^2\,\,\,.
\ee
\end{lemma}
\Proof. By $\om^{(\e)}=r\eta^{(\e)} e_\theta$, %\x{S3ReFRgVoEta},
\x{e63} is identical to
\be \lb{e64}
%\frac{2\pi \abs{\al}}{r_0} (r_0-\e)^3 \le
\norm{r^2\eta^{(\e)}(t)}_{L^1_x}=\norm{r^2\eta^{(\e)}(0)}_{L^1_x} %\le \frac{2\pi \abs{\al}}{r_0} (r_0+\e)^3
%\le M_1.
\le \frac{27\pi}{4}\abs{\kappa}r_0^2\,\,\,.
\ee
The ``inequality'' part of \x{e64} follows from \x{e56} and \x{e61}.
It remains to prove the ``equality'' part, which is actually the conservation of momentum.

Since the initial vorticity field $\om_0^{(\e)}$ in \x{e57} is smooth and compactly supported,
the vorticity field $\om^{(\e)}$ remains Schwartz (smooth and having fast decay in all spatial derivatives) for all the time.
Therefore the momentum can be defined by using the vorticity as
\[
\frac{1}{2} \int_{\R^3} \Bigl (x \times \om^{(\e)}(x,t) \Bigr )\,\,dx\,\,\,,
\]
and moreover, the momentum conserved globally in time, that is
\be \lb{e65}
\frac{1}{2} \int_{\R^3} \Bigl (x \times \om^{(\e)}(x,t) \Bigr )\,\,dx=\frac{1}{2} \int_{\R^3} \Bigl (x \times \om^{(\e)}(x,0) \Bigr )\,\,dx\,,
\hspace{20pt} {\rm for \ all \ t>0}\,,
\ee
which can be checked by the vorticity equations \x{vort}, integration by parts and Schwartz property of the vorticity field $\om^{(\e)}$.

By $\om^{(\e)}=r\eta^{(\e)} e_\theta$, %\x{S3ReFRgVoEta},
\[
\begin{split}
&x \times \om^{(\e)}=x \times r\eta^{(\e)} e_\theta\\
=&(x_1,x_2,x_3) \times \Bigl (-x_2\eta^{(\e)}, x_1\eta^{(\e)},0 \Bigr) \\
=&\Bigl (-x_1x_3\eta^{(\e)}, -x_2x_3\eta^{(\e)},r^2\eta^{(\e)} \Bigr)\,\,. \\
\end{split}
\]
Noting that the first two components are odd in $x_1$ and $x_2$, respectively,
we thus have
\be \lb{e66}
\int_{\R^3}\Bigl (x \times \om^{(\e)}(x,t) \Bigr )\,dx=\biggl (0\,,0\,, \int_{\R^3}r^2\eta^{(\e)}(x,t)\,dx \biggr )\,\,,
\ee
which, combining with \x{e65}, implies
\[
\int_{\R^3}r^2\eta^{(\e)}(x,t)\,dx=\int_{\R^3}r^2\eta^{(\e)}(x,0)\,dx\,, \hspace{20pt} {\rm for \ all \ t>0}\,\,.%=\int_{\R^3}r^2\eta_0^{(\e)}dx
\]
Finally by Lemma \ref{t15}, $\eta^{(\e)}(x,t)$ is nonnegative if $\kappa>0$ (or, nonpositive if $\kappa<0$) for
all points $(x,t) \in \R^3 \times [0, \oo)$ and therefore we can get
\[
\int_{\R^3}\abs{r^2\eta^{(\e)}(x,t)}dx=\int_{\R^3}\abs{r^2\eta^{(\e)}(x,0)}dx\,\,.%=\int_{\R^3}\abs{r^2\eta_0^{(\e)}}dx \le M_1
\]
We get \x{e64} and the lemma is proved.
%which implies the total momentum of the fluid flow is in the $z$-direction.
\qed

\begin{remark} \lb{t17}
The lemma says $\norm{r \om^{(\e)}(t)}_{L^1_x} \lesssim \abs{\kappa} r_0^2$.
\x{e66} implies the total momentum of the fluid flow is in the $z$-direction.
This is due to the special structure of axi-symmetric velocities with no swirl.
\end{remark}

The following lemma claims that the $L^1_x$ norms of $\frac{ \om^{(\e)}}{r}$ are uniformly bounded from above,
which thus gives us the second a-priori estimate.
\begin{lemma} \lb{t18}
For all $t \ge 0$, we have,
\[
\norm{\frac{\om^{(\e)}(t)}{r}}_{L^1_x} \le 3\pi \abs{\kappa}\,\,.
\]
\end{lemma}
\Proof. By $\om^{(\e)}=r\eta^{(\e)} e_\theta$, %\x{S3ReFRgVoEta},
it suffices to prove
\[
\norm{\eta^{(\e)}(t)}_{L^1_x} \le 3\pi \abs{\kappa}, \ \ \ {\rm for \ all} \ t \ge 0\,\,.
\]

We just prove the case of $\kappa>0$. The case of $\kappa<0$ can be proved similarly.
By Lemma \ref{t15}, $\eta^{(\e)} \ge 0$,
Direct calculation shows that
\[
\begin{split}
&\frac{d}{dt} \norm{\eta^{(\e)}(t)}_{L^1_x(\mRTh)}=\frac{d}{dt} \int_\mRTh \eta^{(\e)} \Bigl (x_1,x_2,z,t \Bigr )\,dx_1\,dx_2\,dz \\
=&\int_\mRTh \Bigl (\Delta \eta^{(\e)}- u^{(\e)} \cdot \nabla \eta^{(\e)}+\frac{2}{r}\eta^{(\e)}_{,r} \Bigr )\,dx_1\,dx_2\,dz
=\int_\mRTh \frac{2}{r}\eta^{(\e)}_{,r}\,dx_1\,dx_2\,dz \\
=&4\pi \int_{-\oo}^\oo \int_0^\oo \eta^{(\e)}_{,r} (r,z,t)\,dr\,dz
= -4\pi \int_{-\oo}^\oo \eta^{(\e)} (0,z,t)\,dz \le 0\,\,. \\
\end{split}
\]
Thus $\norm{\eta^{(\e)}(t)}_{L^1_x}$ is decreasing in time.
Combining this with \x{e56}, we get
\[
\norm{\eta^{(\e)}(t)}_{L^1_x} \le \norm{\eta^{(\e)}(0)}_{L^1_x}= \norm{\eta_0^{(\e)}}_{L^1_x} \le 3\pi \abs{\kappa}\,\,.
\]
The lemma is proved. \qed \\

By Nash's method, we will now get uniform estimates of the $L^p_x$ norms of $\frac{\om^{(\e)}}{r}$,
for all $1 \le p \le \oo$, which also serve as our a-priori estimates.
%The ingredients of Nash's method are energy method and Nash's inequality.
%Here we generalize his idea a little bit, considering general $L^p_x$ norms of $\frac{\om^{(\e)}}{r}$ instead of just considering the $L^2_x$ norm.
%Then we get a iteration scheme, by which and by Lemma \ref{t18} we get the desired estimates.
This generalization has been further generalized in \cite{FS86}. %and \cite{K1994}.
The key point in the proof below is that the drift
term $\frac{2}{r}\eta^{(\e)}_{,r}$ %in the equation \x{e62} 
has a good sign.
% and thus can be omitted in the process of proof.
\begin{lemma} \lb{t19}
For every $1 \le p \le \oo$, we have,
\be \lb{e67}
\norm{\frac{\om^{(\e)}(t)}{r}}_{L^p_x} \le C_p t^{-\frac{3}{2}(1-\frac{1}{p})}, \ \ t\in (0,\oo)\,\,,
\ee
where %$C_1$ is from \x{e56}. All
the constants $C_p$ are independent of $\e$.
\end{lemma}
\Proof. %We note that \x{e68} is valid for $p=1$ with $M_1=C_1$ by Lemma \ref{t18}.
%We shall prove, for some constant $C_2$, independent of $\e$, $s$ and $w$,
%\be \lb{decayofGamma}
%\norm{\Gamma^{(\e)}(\cdot,\cdot,\cdot,t;s,w)}_{L^\oo(\mRTh)}\le C_2 t^{-\frac{3}{2}}.
%\ee
%Once we prove \x{decayofGamma}, then by \x{representation}, we have,
%\[
%\begin{split}
%\norm{\eta^{(\e)}(t)}_{L^\oo(\mRTh)} &\le \norm{\Gamma^{(\e)} \Bigl ( \cdot,\cdot,\cdot,t;\sqrt{y_1^2+y_2^2},w \Bigr )}_{L^\oo(\mRTh)}
%\norm{\eta_0^{(\e)}}_{L^1(\mRTh)} \\
%& \le C_1 C_2 t^{-\frac{3}{2}}
%\end{split}
%\]
%Hence \x{e68} is valid for $p=\oo$ with $M_\oo=C_1 C_2$.
%We can then prove \x{e68} for general $0<p<\oo$ by interpolation inequality with
%$M_p=M_1^{1/p}M_\oo^{1-1/p}$. Thus it suffices to prove \x{decayofGamma}.
%To this end, we define the energy
%\[
%E^{(\e)}(t)=\int_\mRTh \abs{\Gamma^{(\e)}(x_1,x_2,z,t;s,w)}^2 dx_1dx_2dz.
%\]
Note that \x{e67} is valid for $p=1$ with $C_1=3\pi\abs{\kappa}$ by Lemma \ref{t18}.
Again by $\om^{(\e)}=r\eta^{(\e)} e_\theta$, %\x{S3ReFRgVoEta},
it suffices to prove
\be \lb{e68}
\norm{\eta^{(\e)}(t)}_{L^p_x} \le C_p t^{-\frac{3}{2}(1-\frac{1}{p})}, \ \ t\in (0,\oo)\,\,.
\ee

Under the spirit of the energy method, for $p=2^n$ with nonnegative integers $n$, we define
\[
E^{(\e)}_p(t) \defnv \norm{\eta^{(\e)}(t)}_{L^p_x}^p=\int_\mRTh \abs{\eta^{(\e)}(x,t)} ^p\, dx\,\,.
\]
For $p=2^n$ with $n \ge 1$, direct calculation yields that
\be \lb{e69}
\begin{split}
-\frac{dE^{(\e)}_p}{dt}=&-\frac{d}{dt}\int_\mRTh \abs{\eta^{(\e)}} ^p\,dx=-\frac{d}{dt}\int_\mRTh \Bigl (\eta^{(\e)} \Bigr )^p\,dx
=-\int_\mRTh p \Bigl (\eta^{(\e)}\Bigr )^{p-1} \eta^{(\e)}_t\,dx \\
=&-\int_\mRTh p \Bigl (\eta^{(\e)}\Bigr )^{p-1} \Bigl (\Delta \eta^{(\e)}+\frac{2}{r} \eta^{(\e)}_{,r}-u^{(\e)} \nabla \eta^{(\e)} \Bigr )\,dx \\
=&-\int_\mRTh \biggl \{p \Bigl [\eta^{(\e)}\Bigr ]^{p-1} \Delta \eta^{(\e)}+\frac{2}{r} \Bigl [(\eta^{(\e)})^p\Bigr ]_{,r}
        - u^{(\e)} \nabla  \Bigl [ (\eta^{(\e)})^p \Bigr ] \biggr \}\,dx \\
=&\int_\mRTh p(p-1) \Bigl [\eta^{(\e)}\Bigr ]^{p-2} \abs{\nabla \eta^{(\e)}}^2\,dx-4\pi\int_{-\oo}^{\oo}\int_{0}^{\oo} \Bigl [(\eta^{(\e)})^p\Bigr ]_{,r}\,dr\,dz \\
=&\int_\mRTh p(p-1) \abs{\Bigl [\eta^{(\e)}\Bigr ]^{\frac{p-2}{2}}\nabla \eta^{(\e)}}^2 dx-4\pi\int_{-\oo}^{\oo} \Bigl [(\eta^{(\e)})^p\Bigr ]^{r=\oo}_{r=0}\,dz \\
=&\int_\mRTh p(p-1) \abs{\frac{2}{p} \nabla  \Bigl [ (\eta^{(\e)})^{\frac{p}{2}} \Bigr ]}^2 dx + 4 \pi \int_{-\oo}^{\oo} \Bigl [(\eta^{(\e)})^p\Bigr ]_{r=0}\,dz \\
\ge &\frac{4(p-1)}{p} \int_\mRTh \abs{\nabla  \Bigl [ (\eta^{(\e)})^{\frac{p}{2}} \Bigr ]}^2 dx\,\,.
\end{split}
\ee
Recall the Nash's inequality \cite[P936]{N58}
\be \lb{e70}
\int_\mRTh \abs{\nabla u}^2 dx \ge M \Bigl ( \int_\mRTh \abs{u}\Bigr )^{-\frac{4}{3}} \Bigl ( \int_\mRTh \abs{u}^2\Bigr )^{\frac{5}{3}}\,\,.
\ee
For $p=2^n$ with $n \ge 1$, by Nash's inequality, we get the following iteration scheme from \x{e69},
\be \lb{e71}
\begin{split}
-\frac{dE^{(\e)}_p}{dt} &\ge \frac{4(p-1)}{p} \int_\mRTh \abs{\nabla  \Bigl [ (\eta^{(\e)})^{\frac{p}{2}} \Bigr ]}^2 dx \\
&\ge \frac{4(p-1)}{p} M \Bigl ( \int_\mRTh \abs{(\eta^{(\e)})^{\frac{p}{2}}}\Bigr )^{-\frac{4}{3}} \Bigl ( \int_\mRTh \abs{(\eta^{(\e)})^{\frac{p}{2}}}^2\Bigr )^{\frac{5}{3}} \\
&=\frac{4(p-1)}{p} M \Bigl ( E^{(\e)}_{p/2}\Bigr )^{-\frac{4}{3}} \Bigl ( E^{(\e)}_{p}\Bigr )^{\frac{5}{3}}\,\,.
\end{split}
\ee
We first prove \x{e68} for $p=2^n$ with nonnegative integers $n$ by induction.
Assume \x{e68} is valid for $q=2^k$ with $k \ge 0$. Let $p=2^{k+1}$. By \x{e71}, we have,
\[
-\frac{dE^{(\e)}_p}{dt} \ge \frac{4(p-1)}{p} M \Bigl ( E^{(\e)}_{q}\Bigr )^{-\frac{4}{3}} \Bigl ( E^{(\e)}_{p}\Bigr )^{\frac{5}{3}}
\ge \frac{4(p-1)}{p} M \Bigl (C^q_q t^{-\frac{3}{2}(q-1)}\Bigr )^{-\frac{4}{3}} \Bigl ( E^{(\e)}_{p}\Bigr )^{\frac{5}{3}}\,\,,
\]
$\Longrightarrow$
\[ \frac{3}{2} \Bigl [ (E^{(\e)}_p)^{-\frac{2}{3}} \Bigr]_t =
-\frac{\frac{dE^{(\e)}_p}{dt}}{\Bigl ( E^{(\e)}_{p}\Bigr )^{\frac{5}{3}}} \ge \frac{4(p-1)}{p} M C_q^{-\frac{4q}{3}}t^{2(q-1)}
=\frac{4(p-1)}{p} M C_q^{-\frac{2p}{3}}t^{p-2}\,\,,
\]
$\Longrightarrow$
\[
(E^{(\e)}_p)^{-\frac{2}{3}}(t) \ge (E^{(\e)}_p)^{-\frac{2}{3}}(t)-(E^{(\e)}_p)^{-\frac{2}{3}}(0)
\ge \frac{8(p-1)}{3p} M C_q^{-\frac{2p}{3}} \int_0^t s^{p-2} ds
= \frac{8M}{3p} C_q^{-\frac{2p}{3}} t^{p-1}\,\,,
\]
$\Longrightarrow$
\[
\norm{\eta^{(\e)}(t)}_{L^p_x}=E^{(\e)}_p(t)^{\frac{1}{p}}
\le \Bigl (\frac{3p}{8M} \Bigr )^{\frac{3}{2p}} C_q t^{-\frac{3}{2}(1-\frac{1}{p})}\,\,.
\]
Hence \x{e68} is valid for $p=2^{k+1}$ with $C_p=\Bigl (\frac{3p}{8M} \Bigr )^{\frac{3}{2p}}C_q$.
In fact, $C_p$ is uniformly bounded from above.
\[
C_p
=\Bigl (\frac{3}{8M} \Bigr )^{\frac{3}{2^{k+2}}}2^{\frac{3(k+1)}{2^{k+2}}}C_{2^k}
\le \Bigl (\frac{3}{8M} \Bigr )^{\sum \frac{3}{2^{i+2}}}2^{\sum \frac{3(i+1)}{2^{i+2}}}C_1 =: C_\oo\,\,.
\]
$\Longrightarrow$
\[ \norm{\eta^{(\e)}(t)}_{L^\oo_x} \le C_\oo t^{-\frac{3}{2}}\,\,. \]
For other $p$, we can prove \x{e68} by interpolation. Therefore the lemma is proved. \qed

\begin{remark} \lb{t20}
From the proof of Lemma \ref{t19}, we see the constants $C_p$ in \x{e67}
linearly depends on $C_1=3\pi \abs{\kappa}$. In particular, %by \x{e56}, we have
\be \lb{e72}
\begin{split}
&C_\oo=\Bigl (\frac{3}{8M} \Bigr )^{\sum \frac{3}{2^{i+2}}}2^{\sum \frac{3(i+1)}{2^{i+2}}}C_1 \lesssim \abs{\kappa}\,\,, \\
&\norm{\frac{\om^{(\e)}(t)}{r}}_{L^\oo_x} \le C_\oo t^{-\frac{3}{2}} \lesssim \abs{\kappa} t^{-\frac{3}{2}}\,\,,
\end{split}
\ee
which gives us the third a-priori estimate, where $M$ is the absolute constant in Nash's inequality \x{e70}.
\end{remark}
%\begin{remark} \lb{NashtoGMO}
%Nash's method also works for the problem of 2d NSE with finite measures as initial vorticity \cite{GMO},
%which has been observed by Kato \cite{K1994}.
%Let 2D initial vorticity $\om_0$ be an arbitrary finite measure. Let $\om^{(\e)}$
%be the vorticity field related to the regularized initial data $\om_0^{(\e)}=\phi^{(\e)}*\om_0$.
%After applying the above argument in Lemma \ref{t19} to $\om^{(\e)}$, one shall get,
%\[
%\norm{\om^{(\e)}(t)}_{L^p_x(\mRT)} \le M_p t^{-(1-\frac{1}{p})}.
%\]
%\end{remark}
\begin{remark} \lb{t21}
If the fluid is inviscid, then $\eta^{(\e)}$ satisfies
\be \lb{e73}
\eta_t+ u^{(\e)} \cdot \nabla \eta=0, \hspace{10pt} {\rm in} \ \ \mRTh \times (0, \oo)\,\,.
\ee
Since $\eta^{(\e)}$ is conserved along particle trajectories, $\eta^{(\e)}$ keeps its sign in later time.
We still have the uniform estimates of the $L^1_x$ norms:
\[
\norm{\eta^{(\e)}(t)}_{L^1_x} = \norm{\eta^{(\e)}(0)}_{L^1_x}= \norm{\eta_0^{(\e)}}_{L^1_x} \le 3\pi\abs{\kappa}\,\,.
\]
However, the argument in Lemma \ref{t19} yields: for any $1<p \le \oo$,
\[
\norm{\eta^{(\e)}(t)}_{L^p_x}=\norm{\eta^{(\e)}(0)}_{L^p_x}=\norm{\eta_0^{(\e)}}_{L^p_x}\,\,,
\]
which will blow up as $\e$ goes to 0. Therefore we lose uniform controls of the $L^p_x$ norms in the inviscid case.
\end{remark}

We now use the weighted inequalities of the previous section and
the three a-priori estimates from Lemma \ref{t16}, Lemma \ref{t18} and Remark \ref{t20}
to get further estimates on
vorticity, the gradient of velocity, velocity and pressure.
\begin{lemma} \lb{t22} For $0<t<\oo$, we have the following estimates:\\
{\rm i)} for any $1 \le p \le 2$
\be \lb{e74}
\norm{\om^{(\e)}(t)}_{L^p_x} \lesssim \abs{\kappa}r_0 t^{-\frac{3}{2}\bigl( 1-\frac{1}{p}\bigr)},
\ee
{\rm ii)} for any $1 < p \le 2$
\be \lb{e75}
\norm{\nabla u^{(\e)}(t)}_{L^p_x} \lesssim \abs{\kappa}r_0 t^{-\frac{3}{2}\bigl( 1-\frac{1}{p}\bigr)},
\ee
{\rm iii)} for any $\frac{3}{2} < q \le 6$
\be \lb{e76}
\norm{u^{(\e)}(t)}_{L^q_x} \lesssim \abs{\kappa}r_0t^{-\bigl( 1-\frac{3}{2q}\bigr)},
\ee
{\rm iv)} for any $1 < q \le 3$
\be \lb{e77}
\norm{p^{(\e)}(t)}_{L^q_x} \lesssim \abs{\kappa}^2 r_0^2 t^{-\bigl( 2-\frac{3}{2q}\bigr)},
\ee
{\rm v)}
\be \lb{e78}
\norm{u^{(\e)}(t)}_{L^\oo_x} \lesssim \abs{\kappa} r_0^\frac{1}{2} t^{-\frac{3}{4}}\,\,.
\ee
\end{lemma}
\Proof.

i). By Proposition \ref{t2}, for any $1 \le p \le 2$, we have
\be \lb{e79}
\norm{\om^{(\e)}(t)}_{L^p_x} \le \norm{r \om^{(\e)}(t)}_{L^1_x}^{\frac{1}{2}} \norm{\frac{\om^{(\e)}(t)}{r}}_{L^1_x}^{\frac{1}{p}-\frac{1}{2}} \norm{\frac{ \om^{(\e)}(t)}{r}}_{L^\oo_x}^{1-\frac{1}{p}}\,\,.
\ee
Then \x{e74} is an easy consequence of \x{e79}, Lemma \ref{t16}, Lemma \ref{t18} and \x{e72} in Remark \ref{t20}.

{\rm ii)} By div\,$u^{(\e)}$=0, curl\,$u^{(\e)}$=$\om^{(\e)}=(\om^{(\e)}_1,\om^{(\e)}_2,\om^{(\e)}_3)$ and Fourier transform,
one can get
\[
\nabla u^{(\e)} \\
=\begin{bmatrix}
R_1R_2\om^{(\e)}_3-R_1R_3\om^{(\e)}_2&R_2R_2\om^{(\e)}_3-R_2R_3\om^{(\e)}_2&R_2R_3\om^{(\e)}_3-R_3R_3\om^{(\e)}_2 \\
R_1R_3\om^{(\e)}_1-R_1R_1\om^{(\e)}_3&R_2R_3\om^{(\e)}_1-R_1R_2\om^{(\e)}_3&R_3R_3\om^{(\e)}_1-R_1R_3\om^{(\e)}_3 \\
R_1R_1\om^{(\e)}_2-R_1R_2\om^{(\e)}_1&R_1R_2\om^{(\e)}_2-R_2R_2\om^{(\e)}_1&R_1R_3\om^{(\e)}_2-R_2R_3\om^{(\e)}_1
\end{bmatrix}\,\,,
\]
where $R_j$, $j=1, 2, 3$ are the classical Riesz transformations, which are well-defined and continuous
on $L^p_x(\R^3)$ for all $1<p<\oo$, see for instance \cite{Stein93,Torchinsky86}. Therefore
\[ %\lb{vorticitycontrolvelocitygradient}
\norm{\nabla u^{(\e)}(t)}_{L^p_x} \lesssim \norm{\om^{(\e)}(t)}_{L^p_x}\,\,,
\]
which, combining with \x{e74},
implies \x{e75}.

{\rm iii)} By Corollary \ref{t3}, for any $\frac{3}{2}<q \le 6$,
\be \lb{e80}
\norm{u^{(\e)}(t)}_{L^q_x} \lesssim \norm{r \om^{(\e)}(t)}_{L^1_x}^{\frac{1}{2}}\norm{\frac{\om^{(\e)}(t)}{r}}_{L^1_x}^{\frac{1}{q}-\frac{1}{6}}\norm{\frac{ \om^{(\e)}(t)}{r}}_{L^\oo_x}^{\frac{2}{3}-\frac{1}{q}}\,\,.
\ee
Then \x{e76} is an easy consequence of \x{e80}, Lemma \ref{t16}, Lemma \ref{t18} and \x{e72}.

{\rm iv)} Recall the pressure $p^{(\e)}$ and the velocity $u^{(\e)}=\Bigl (u^{(\e)}_1, u^{(\e)}_2, u^{(\e)}_3 \Bigr )$
satisfy the following equation (which can be easily obtained from Navier-Stokes equations and divergence-free condition div\,$u^{(\e)}$=0):
\be \lb{e81}
\Delta p^{(\e)}= -\partial_j \partial_k (u^{(\e)}_ju^{(\e)}_k)\,\,.
\ee
Then by \x{e76},
we can use the Riesz transformation $R_j$ to solve
\x{e81} to get
\[
p^{(\e)}=R_jR_k(u^{(\e)}_ju^{(\e)}_k)\,\,.
\]
Hence
\[ %\lb{velocitycontrolpressure}
\norm{p^{(\e)}(t)}_{L^q_x} \lesssim \norm{u^{(\e)}(t)}_{L^{2q}_x}^2\,\,,
\]
which, combining with \x{e76}, implies \x{e77}.

{\rm v)} By Proposition \ref{t10} and \ref{t12}, %for any $\kappa \in \Bigl (\frac{1}{2}, \frac{2}{3} \Bigr]$,
\be \lb{e82}
\begin{split}
\norm{u^{(\e)}(t)}_{L^\oo_x}
\lesssim \ &\norm{r \om^{(\e)}(t)}_{L^1_x}^{\frac{1}{4}} \norm{\frac{\om^{(\e)}(t)}{r}}^{\frac{1}{4}}_{L^1_x} \norm{\frac{\om^{(\e)}(t)}{r}}^{\frac{1}{2}}_{L^\oo_x}\,\,.
\end{split}
\ee
Then \x{e78} is an easy consequence of \x{e82}, Lemma \ref{t16}, Lemma \ref{t18} and \x{e72}.
\qed

By Lemma \ref{t22} and the subcritical theory of Navier-Stokes equations,
we can control the spatial and time derivatives of the velocity and pressure of any order pointwise.
\begin{lemma} \lb{t23} For any $k, h \ge 0$ and for any $0<s<T$, we have the following pointwise estimate
\[ \norm{\nabla^k_x\nabla^h_t u^{(\e)}}_{C^0_{x,t}(\R^3 \times [s,T])} \le C, \hspace{20pt}
\norm{\nabla^k_x\nabla^h_t p^{(\e)}}_{C^0_{x,t}(\R^3 \times [s,T])} \le C\,\,,
\]
where $C$ is independent of $\e$ and depends only on $k,h, s,T,\abs{\kappa},r_0$.
\end{lemma}
\Proof. This lemma is a consequence of the subcritical well-posedness theory of Navier-Stokes equations. Fix $0<s<T$. By \x{e76},
we have the following subcritical estimate
\be \lb{e83}
\norm{u^{(\e)} (t)}_{L^6_x} \lesssim \abs{\kappa}r_0t^{-\frac{3}{4}}\,\,,
\ee
since $L^6_x(\R^3)$ is a subcritical space for Navier-Stokes equations with respect to the scaling
\[
u(x, t) \longmapsto \la u(\la x, \la^2 t), \ p(x, t)\longmapsto \la^2 p(\la x, \la^2 t)\,\,.
\]
By the standard subcritical theory, see for instance \cite{K84, GMO}, there exists a local-in-time unique solution $v^{(\e)}$ for
Navier-Stokes equations with $u^{(\e)} \Bigl (\frac{s}{2} \Bigr )$ as initial velocity in the space $C\Bigl (\bigl [\frac{s}{2}, T^* \bigr ), L^6_x(\R^3) \Bigr )$ for some $\frac{s}{2}<T^* \le \oo$. $v^{(\e)}$ coincides with $ u^{(\e)}$ on
the time interval $[\frac{s}{2},T^*)$ by weak-strong uniqueness.
The decay property \x{e83} implies $T^*=\oo$.
Hence $u^{(\e)}=v^{(\e)}$ for all $t \in [\frac{s}{2}, \oo)$. Again by the subcritical theory,
$u^{(\e)}$ satisfies
\be \lb{e84}
\norm{\nabla^k_x\nabla^h_t u^{(\e)}}_{L^\oo_tL^6_x(\R^3 \times [s,T])} \le C\,\,,
\ee
where $C$ depends only $k, h, s, T, \norm{u^{(\e)} \Bigl (\frac{s}{2} \Bigr )}_{L^6_x}$.
Then by Sobolev embedding, we prove the first estimate. The second estimate is a consequence of \x{e84} and \x{e81}.
\qed

The estimate \x{e76} in Lemma \ref{t22} imply
the set $\big \{u^{(\e)} \bigr \}_{0<\e<\frac{r_0}{2}}$ has weak compactness in Lebesgue spaces.
To show the strong convergence of $\big \{ u^{(\e)} \bigr \}_{0<\e<\frac{r_0}{2}}$,
we need to establish certain uniform weak continuity of $u^{(\e)}$
as functions of time $t$. To this end, we follow the standard method,
see \cite{Cottet1986, Lin98,T77}. Let $H^{-2}_x(\R^3)$ be the dual space of $H^2_x(\R^3)$.
\begin{lemma} \lb{t24}
Let $0<T<\oo$. Then we have
\be
\norm{\frac{\partial u^{(\e)}}{\partial t}}_{L^{\frac{5}{4}}_t \bigl ( 0, T; H^{-2}_x(\R^3) \bigr )} \le C\,\,,
\ee
where the constant $C$ is independent of $\e$ and depends on $T$.
\end{lemma}
\Proof. Let $\phi \in H^2_x(\R^3)$. By Navier-Stokes equations and Lemma \ref{t22}, we have
\[
\begin{split}
&\abs{\Bigl (\frac{\partial u^{(\e)}}{\partial t}, \phi \Bigr )} \\
=&\abs{\Bigl( -{\rm div}(u^{(\e)} \otimes u^{(\e)})-\nabla p^{(\e)}+\Delta u^{(\e)}, \phi \Bigr )} \\
\le &\abs{(u^{(\e)} \otimes u^{(\e)}, \nabla \phi)}+\abs{(p^{(\e)}, {\rm div} \phi)}+\abs{(u^{(\e)}, \Delta \phi)} \\
\lesssim &\norm{u^{(\e)}(t)}_{L^{p_1}_x} \norm{u^{(\e)}(t)}_{L^{p_2}_x} \norm{\nabla \phi }_{L^{p_3}_x}
+\norm{p^{(\e)}(t)}_{L^{q_1}_x}  \norm{\nabla \phi }_{L^{q_2}_x}+\norm{u^{(\e)}(t)}_{L^2_x} \norm{\Delta \phi }_{L^2_x} \\
\lesssim & \abs{\kappa}^2 r_0^2 t^{-2+\frac{3}{2p_1}+ \frac{3}{2p_2}} \norm{\nabla \phi}_{L^{p_3}_x}+ \abs{\kappa}^2 r_0^2 t^{-2+\frac{3}{2q_1}} \norm{\nabla  \phi}_{L^{q_2}_x}
+\abs{\kappa} r_0t^{-\frac{1}{4}} \norm{\phi}_{H^2_x}\,\,,
\end{split}
\]
where
\[
\frac{1}{p_1}+\frac{1}{p_2}+\frac{1}{p_3}=1, \ \frac{3}{2} < p_1, p_2 \le 6,\ \frac{1}{q_1}+\frac{1}{q_2}=1, \ 1<q_1 \le 3\,\,.
\]
One can take, for example,
\[
p_1=p_2=\frac{12}{5},\ p_3=6,\ q_1=\frac{6}{5},\ q_2=6\,\,.
\]
Then by Sobolev embedding, we have for $0<t \le T$,
\[
\begin{split}
&\abs{\Bigl (\frac{\partial u^{(\e)}}{\partial t}, \phi \Bigr )} \\
\lesssim &\abs{\kappa}^2 r_0^2 t^{-\frac{3}{4}} \norm{\nabla \phi}_{L^{6}_x}
+\abs{\kappa}^2 r_0^2 t^{-\frac{3}{4}} \norm{\nabla \phi}_{L^{6}_x}+
\abs{\kappa} r_0 t^{-\frac{1}{4}} \norm{\phi}_{H^2_x} \\
\lesssim &  \Bigl (\abs{\kappa}^2 r_0^2 t^{-\frac{3}{4}}+\abs{\kappa} r_0 t^{-\frac{1}{4}}\Bigr ) \norm{\phi}_{H^2_x}\,\,.
\end{split}
\]
Hence
\[
\norm{\frac{\partial u^{(\e)}}{\partial t}(t)}_{H^{-2}_x} \lesssim \abs{\kappa}^2 r_0^2 t^{-\frac{3}{4}}+\abs{\kappa} r_0 t^{-\frac{1}{4}}\,\,.
\]
Finally integrating with respect to time from $(0,T)$ yields the desired result.
\qed
\begin{lemma} \lb{t25}
For any $0<T<\oo$, $\Bigl \{ u^{(\e)} \Bigr \}_{0<\e<\frac{r_0}{2}}$ is precompact in $L^{\frac{8}{5}}_t \bigl (0,T; L^2_{x,{\rm loc}}(\R^3) \bigr)$.
\end{lemma}
\Proof. By Lemma \ref{t22}, one has
\[
\norm{u^{(\e)}(t)}_{L^{\frac{8}{5}}_x} \lesssim  \abs{\kappa} r_0 t^{-\frac{1}{16}},
\hspace{20pt}  \norm{\nabla u^{(\e)}(t)}_{L^{\frac{8}{5}}_x} \lesssim  \abs{\kappa} r_0 t^{-\frac{9}{16}}\,\,,
\]
which implies
\be \lb{e85}
u^{(\e)} \in L^{\frac{8}{5}}_t \bigl (0,T; W^{1,\frac{8}{5}}_x(\R^3) \bigr)\,\,.
\ee
Then \x{e85}, Lemma \ref{t24} and Theorem 2.1 of \cite[Chap. III]{T77}
imply the desired result.
\qed

\section{Proof of Theorem \ref{t1}}
\x{e76} implies
\be \lb{e86}
\norm{u^{(\e)}(t)}_{L^2_x} \lesssim \abs{\kappa} r_0 t^{-\frac{1}{4}}\,\,,
\ee
which in turn implies
\be \lb{e87}
\set{u^{(\e)}} \ {\rm is \ a \ bounded \ set \ in} \ L^{\frac{8}{3}}_t L^2_x \bigl (\R^3 \times (0,T) \bigr), \hspace{20pt} {\rm for \ any} \ 0<T<\oo\,\,.
\ee
Arzela-Ascoli's Theorem, Lemma \ref{t23}, Lemma \ref{t25} and \x{e87}
allow us to extract a subsequece of $\set{u^{(\e)}, p^{(\e)}}$,
still denoted as $\set{u^{(\e)}, p^{(\e)}}$ such that for a smooth vector field $u$ and a smooth scalar
function $p$, for any nonnegative integers $k,h$ and for any $0<T<\oo$, we have
\be \lb{e88}
\begin{split}
&u^{(\e)} \to u \ \ \ {\rm in} \ L^{\frac{8}{5}}_t \bigl (0,T; L^2_{x, {\rm loc}}(\R^3) \bigr)\,\,, \\
&u^{(\e)} \rightharpoonup u \ \ \ {\rm in} \ L^{\frac{8}{3}}_t L^2_x \bigl (\R^3 \times (0,T) \bigr)\,\,,
\end{split}
\ee
and
\be \lb{e89}
\begin{split}
&\nabla^k_x \nabla^h_t u^{(\e)} \rightrightarrows \nabla^k_x \nabla^h_t u \ \ \ {\rm locally \ in} \R^3 \times (0,\oo)\,\,, \\
%{\rm in} \ L^\oo_{{\rm loc}}(\R^3 \times (0,\oo)), \\
&\nabla^k_x \nabla^h_t p^{(\e)} \rightrightarrows \nabla^k_x \nabla^h_t p \ \ \ {\rm locally \ in} \R^3 \times (0,\oo)\,\,,
%{\rm in} \ L^\oo_{{\rm loc}}(\R^3 \times (0,\oo)),
\end{split}
\ee
which imply the limit $\bigl (u, p \bigr )$ is a global-in-time smooth
solution of Navier-Stokes equations in $\R^3 \times (0,\oo)$ and $u$ is axi-symmetric with no swirl.

We prove the initial condition \x{e1}. Take a $\vf \in C^\oo_0 \bigl (\R^3;\R^3 \bigr )$ with $B_R(0)$ as its support.
By Navier-Stokes equations, we have
\be \lb{e90}
\begin{split}
&\int_0^T\int_{\R^3} \biggr \{
\bigl (u^{(\e)} \otimes u^{(\e)} \bigr) \cdot \nabla {\rm curl} \vf+
u^{(\e)} \cdot \Delta {\rm curl} \vf
\biggr \}\,dx\,dt \\
&=\int_{\R^3} \om^{(\e)}(x, T) \cdot \vf(x)\,dx-\int_{\R^3} \om^{(\e)}_0(x) \cdot \vf(x)\,dx\,\,.
\end{split}
\ee
We claim that we are able to pass to the limit in \x{e90} to get
\be \lb{e91}
\begin{split}
&\int_0^T\int_{\R^3} \biggr \{\bigl (u \otimes u \bigr) \cdot \nabla {\rm curl} \vf+
u \cdot \Delta {\rm curl}\vf
\biggr \}\,dx\,dt \\
&=\int_{\R^3} \om(x,T) \cdot \vf(x)\,\,dx-\int_{\R^3} \ourring \cdot \vf \,\, dx\,\,.
\end{split}
\ee
To this end, it suffices to check the nonlinear term in \x{e90} and \x{e91}.
By \x{e87} and \x{e88}, we have
\[
\begin{split}
&\abs{\int_0^T\int_{\R^3} \bigl (u^{(\e)} \otimes u^{(\e)} \bigr) \cdot \nabla {\rm curl} \vf \,dx\,dt
-\int_0^T\int_{\R^3}  \bigl ( u \otimes u \bigr) \cdot \nabla {\rm curl} \vf\,dx\,dt} \\
\lesssim &\norm{u^{(\e)}-u}_{L^{\frac{8}{5}}_t L^2_x \bigl (B_R(0) \times (0,T) \bigr)}
\biggl ( \norm{u^{(\e)}}_{L^{\frac{8}{3}}_t L^2_x \bigl (\R^3 \times (0,T) \bigr)}+\norm{u}_{L^{\frac{8}{3}}_t L^2_x \bigl (\R^3 \times (0,T) \bigr)} \biggr )
\norm{\nabla {\rm curl} \vf}_{L^\oo_x},
\end{split}
\]
which goes to $0$ as $\e \to 0$. Thus \x{e91} is obtained.

Fatou's lemma and \x{e86} imply
\be
\norm{u(t)}_{L^2_x} \lesssim \abs{\kappa} r_0 t^{-\frac{1}{4}}\,\,.
\ee
Hence in view of \x{e91}
\be
\begin{split}
&\abs{\int_{\R^3} \om(x,T) \cdot \vf(x)\,dx-\int_{\R^3} \ourring \cdot \vf\, dx} \\
=&\abs{\int_0^T\int_{\R^3} \biggr \{\bigl (u \otimes u \bigr) \cdot \nabla {\rm curl} \vf+
u \cdot \Delta {\rm curl} \vf
\biggr \}\,dx\,dt} \\
\lesssim & \int_0^T \abs{\kappa}^2 r_0^2 t^{-\frac{1}{2}}\,dt+\int_0^T \abs{\kappa} r_0 t^{-\frac{1}{4}}\,dt \lesssim \abs{\kappa}^2 r_0^2 T^{\frac{1}{2}}+\abs{\kappa} r_0 T^{\frac{3}{4}},
\end{split}
\ee
which implies \x{e1}. Theorem \ref{t1} is thus proved.

\begin{remark}
Theorem \ref{t1} is also true if we replace the initial condition by
finite many vortex rings
\be \lb{e92}
\om(\cdot,0) = \sum_{i=1}^n \kappa_i \delta_{r_i,z_i} e_\theta\,\,,
\ee
where all $\kappa_i>0$ (or, all $\kappa_i <0$), or more generally, by
\be \lb{e93}
\om(\cdot,0) = \mu e_\theta\,\,,
\ee
where $\mu$ is a positive or negative finite measure with a compact support in the $rz$-plane.
Without any modification,
the preceding proof for single vortex ring also works for the cases of \x{e92}
and \x{e93}.
\end{remark}

\end{document}